\documentclass[12pt,reqno]{amsart}
\usepackage{amsfonts,amsmath,amssymb,amsthm}
\usepackage{latexsym}
\usepackage{eucal}
\usepackage[ansinew]{inputenc}
\usepackage[american]{babel}
\usepackage{amsfonts}
\usepackage{amssymb}
\usepackage[dvips]{graphicx}

\numberwithin{equation}{section}

\newtheorem{teo}{Theorem}[section]
\newtheorem{lema}{Lemma}[section]
\newtheorem{prop}{Proposition}[section]
\newtheorem{defi}{Definition}[section]

\newtheorem{obs}{Remark}[section]
\newtheorem{coro}{Corollary}[section]

\begin{document}

\title[(In) stability of Solitary Waves for fKdV and fBBM]
{Stability properties of solitary waves for  fractional KdV and BBM equations}

\author[J. Angulo]{Jaime Angulo Pava}

\maketitle

\centerline{Department of Mathematics,
IME-USP}
 \centerline{Rua do Mat\~ao 1010, Cidade Universit\'aria, CEP 05508-090,
 S\~ao Paulo, SP, Brazil.}
 \centerline{\it angulo@ime.usp.br}
\begin{abstract}
This paper sheds new light on the stability properties of solitary wave  solutions associated with models of Korteweg-de Vries  and Benjamin\&Bona\&Mahoney type, when the dispersion is very lower. Via an approach of compactness,  analyticity and asymptotic perturbation theory, we establish sufficient conditions for the
existence of exponentially growing solutions to the linearized
problem and so a criterium of  linear instability of solitary waves is obtained for both models. Moreover, the nonlinear stability and linear instability of the ground states solutions for both models is obtained for some specific regimen of parameters. Via a Lyapunov strategy and a variational analysis we obtain the stability of the blow-up 
of solitary waves for the critical fractional KdV equation. 

The  arguments presented in this investigation
has prospects for the study of the instability of
traveling waves solutions  of other nonlinear evolution equations.
\end{abstract}




\thanks{{\it Date}: January 20, 2017}


\section{Introduction}

This paper provides a detailed study of various stability issues associated to the  dynamic of solitary waves solutions for the  so-called  fractional  Korteweg-de Vries equation (henceforth fKdV equation)
\begin{equation}\label{0gfkdv}
u_t+u^pu_x-D^\alpha u_x=0,\qquad p\in \mathbb N,
\end{equation}
where  $u=u(x,t)$, $x,t\in \mathbb{R}$, represents a real
valued function, and $D^\alpha$ is defined via Fourier transform by
$$
\widehat{D^\alpha f}(\xi)=|\xi|^\alpha \widehat{f}(\xi), \qquad \alpha\in (0,1).
$$

The importance of  study of this model for any $\alpha>0$ is due to its physical relevance and its own mathematical interest. We recall that the model (\ref{0gfkdv}) contains two famous family of equations, the generalized Korteweg-de Vries for $\alpha=2$  (gKdV henceforth),  and the generalized Benjamin-Ono equation for  $\alpha=1$ (gBO henceforth), and  in this case $D$ can be write as $D=\mathcal H\partial_x$, where $\mathcal H$  denotes  the Hilbert transform and which may be defined by
$\widehat{\mathcal Hf}(\xi)=-i sgn (\xi) \widehat{f}(\xi)$.

For $\alpha\geqq 1$, studies on the Cauchy problem, blow-up issues, large-time asymptotic behavior of solutions, the stability of solitary waves solutions, breathers solutions and multi-solitons solutions (as well as periodic traveling  wave solutions) have been  the focus  of deep research in the past years via a rich variety of techniques, see by instance, Albert \cite{Al1}, Albert\&Bona \cite{AlB}, Albert\&Bona\&Saut \cite{ABS}, Alejo\&Mu$\bar{n}$oz \cite{AM}, Angulo \cite{angulo4}, Angulo\&Bona\&Scialom \cite{AnBS}, Benjamin \cite{be1}-\cite{Be2}, Bona \cite{Bo1}, Bona\&Souganidis\&Strauss \cite{BSS}, Bona\&Saut \cite{BSa}, Iorio \cite{I1}, Grillakis\&Shatah\&Strauss \cite{grillakis1}-\cite{grillakis2}, Kenig\&Ponce\&Vega \cite{KPV1}, Lopes \cite{lopes2}, Martel\&Merle \cite{MM1}-\cite{MM2}-\cite{MM3}, Martel\&Pilod \cite{MP}, Mu$\bar{n}$oz \cite{Mu1}, Weinstein \cite{W1}-\cite{W2}.

The case $\alpha\in (0,1)$ has been the focus of many recent  studies.  The Cauchy problem, the existence of solitary wave solutions, the stability properties of the ground states and numerical simulations  have been addressed by  Linares\&Pilod\&Saut \cite{fds}-\cite{fds2}, Frank\&Lenzmann \cite{FL} and Klein\&Saut \cite{KS}. 

One of the objectives of this paper is to extend the theory  of Vock\&Hunziker in \cite{VH} about the stability of Schr\"odinger eigenvalue problems to the  study of linear instability of solitary waves solutions for the fKdV equation with a ``lower dispersion'' (see Theorems \ref{teo2}-\ref{teo3} below).  In particular, we recover the linear instability results in \cite{KSt} for the ground state solutions (see Definition \ref{defig} below) of (\ref{0gfkdv}) with $p=1$ and $\alpha \in (\frac 13, \frac12)$. For completeness of the  exposition, we show in an unified way the nonlinear stability results for the ground state solutions  with $p<2\alpha$ and $\alpha \in (\frac12, 2)$ (see Theorem \ref{nonlinear} below). 

The case $\alpha=\frac12$ and $p=1$ in (\ref{0gfkdv}), so-called the {\it critical case} for the fKdV model, it remains open for a stability analysis of solitary waves. Indeed, in this case,  the recently numerical simulations in Klein\&Saut \cite{KS} suggest  the existence of blow up of solutions for initial data close to the solitary waves and proving such result seems to be out of reach. Here, we will show that for this critical case a kind of ``{\it stability of the blow-up}'' near to the possible unstable ground state solutions happens  and we checked  one of the conjectures emerging of the numerical findings in \cite{KS} (see Theorem \ref{t4.9} below).

Our approach of linear instability for the solitary wave solutions of (\ref{0gfkdv}) is extended to the following generalized fractional Korteweg-de Vries models (gfKdV henceforth)
\begin{equation}\label{gfkdv}
u_t+(f(u))_x-(\mathcal{M}u)_x=0,
\end{equation}
where $\mathcal{M}$ is a differential or pseudo-differential operator  defined as a Fourier
multiplier operator
\begin{equation}\label{opera}
\widehat{\mathcal{M}g}(\xi)=\beta (\xi)\widehat{g}(\xi),\;\;\;\; \xi\in \mathbb R,
\end{equation}
and, $f$ is assumed to be a smooth nonlinear function.
The symbol $\beta$ of $\mathcal{M}$ (representing the lower dispersion
effects) is assumed to be continuous, locally bounded, even
function on $\mathbb R$, satisfying the conditions
$$
 a_1|\xi|^{\gamma} \leq \beta(\xi)\leq a_2
(1+|\xi|)^{\alpha},
$$
for   $ |\xi|\geq b_0$, $0\leq \gamma\leqq \alpha<1$,
with $\beta(\xi)>b$, for all $\xi \in \mathbb R$ and $a_i>0$, $i=1,2$. In this point of the analysis, we extend the linear instability results in Lin \cite{lin} for the models (\ref{gfkdv})
with a growth  of the symbol of $\mathcal M$ determined by $\alpha \in (0,1)$. 

 We note, that in various models of fluid dynamics and mathematical physics the symbol $\beta$ in (\ref{opera}) is not necessarily polynomial, such as in the case of the Whitham equation for describing water waves in the small amplitude and long wave regime when surface tension is included (\cite{La}, \cite{LS}, \cite{Wh}) 
\begin{equation}
\beta(\xi)=(1+\gamma |\xi|^2)^{1/2}\Big(\frac{tanh \xi}{\xi}\Big)^{1/2},
\end{equation}
where $\gamma \geqq 0$  measures the surface tension effects. Here $\beta$ satisfies 
$$
\frac12|\xi|^{1/2}\leqq \beta (\xi)\leqq 2|\xi|^{1/2}, \;\;\;\text{for}\;\; |\xi|\;\; \text{large}.
$$

The analysis established above for the fKdV equation (\ref{0gfkdv}) was also extended to the fractional BBM equation (fBBM henceforth)
\begin{equation}\label{fBBM0}
u_t+u_x + \partial_x(u^2)+D^\alpha u_t=0,
\end{equation}
for $\alpha\in (\frac13, 1)$. In this case we show that the ground state solutions associated to the fBBM equation are linearly unstable for $\alpha\in (\frac13, \frac12)$ with a wave-speed not so large. For $\alpha\in (\frac13, 1)$ we have also nonlinear stability in the case of the  wave-speed  in general to be large (see Theorems \ref{nonlinearfbbm1}-\ref{linsfbbm1} below).

Before of establishing  more precisely our results, we will make a brief summary of some known results for the fKdV model,  $\alpha\in (0, 1)$, it which will be useful in our exposition. We start initially with some basic information about the existence of solitary waves solutions for this model. A solitary wave solution for (\ref{0gfkdv}) is a solution of the form $u(x,t)=\varphi_c(x-ct)$ with
$$
\lim_{|\xi|\to \infty} \varphi_c(\xi)=0,
$$
which (if they exist !!!) it will represent a ``perfect'' balance between the lower dispersion and the effects of the nonlinearity. For $\varphi\equiv \varphi_c$ belongs to the space $H^{\alpha/2}(\mathbb R)\cap L^{p+2}(\mathbb R)$ we have that 
\begin{equation}\label{equa}
D^\alpha \varphi +c \varphi- \frac{1}{p+1} \varphi^{p+1}=0.
\end{equation}
The existence of solutions for (\ref{equa}), with the later specify regularity conditions, it can be deduced from the Concentration-Compactness Method for any $c>0$ and  $p\in (1, \frac{2\alpha}{1-\alpha})$ (see Weinstein \cite{W2} and Arnesen \cite{Arn}). We can also to see  (by the so-called Pohozaev identities) that the pseudo-differential equation satisfied by the profile $ \varphi$ does not admit any non-trivial solutions for the following cases:
\begin{enumerate}
\item[(1)] for $\alpha\geqq 1$ and $c<0$ (without restrictions on the power $p$),
\item[(2)] for $\alpha\in (0,1)$, $c>0$ and $\alpha\leqq \frac{p}{p+2}$.
\end{enumerate}
 For completeness of the exposition, we will establish the item (2) above (the item (1) is very well known). Indeed, for $\alpha\in (0,1)$ we have from Lemma B.2 in Frank\&Lenzmann \cite{FL}  that $\varphi\in H^{\alpha +1}(\mathbb R)$. Then by Plancherel Theorem, the following energy identity is immediate
\begin{equation}\label{ener}
\int_{\mathbb R} |D^{\alpha/2} \varphi|^2 dx + c \int_{\mathbb R} \varphi^2dx-\frac{1}{p+1} \int_{\mathbb R}\varphi^{p+2}dx=0.
\end{equation}
Next, since $\varphi'$ makes sense, we  have by Plancherel  and integration by parts that
\begin{equation}\label{poha}
\int_{\mathbb R} x \varphi'D^{\alpha} \varphi dx =\frac{\alpha-1}{2}\int_{\mathbb R} |D^{\alpha/2} \varphi|^2 dx.
\end{equation}
Thus, from (\ref{ener})-(\ref{poha}) follows
\begin{equation}\label{aprior}
(\alpha(p+2)-p)\int_{\mathbb R} |D^{\alpha/2} \varphi|^2 dx=c p\int_{\mathbb R} |\varphi|^2 dx,
\end{equation}
proving that no finite energy solitary waves exist when $c>0$ and $\alpha\leqq \frac{p}{p+2}$ hold. 

The following definition will be useful in our study (see Frank\&Lenzmann \cite{FL}). 

\begin{defi}\label{defig} Let $Q\in H^{\alpha/2}(\mathbb R)$ be an even and positive solution of 
\begin{equation}\label{dground}
D^\alpha Q+Q-Q^{p+1}=0\qquad in \;\;\;\mathbb R.
\end{equation}
If $Q$ solves the minimization problem
\begin{equation}\label{minimi}
J^{\alpha, p}(Q)=inf\{J^{\alpha, p}(v): v\in H^{\alpha/2}(\mathbb R)-\{0\}\}
\end{equation}
where $J^{\alpha, p}$ is the `Weinstein' functional
\begin{equation}\label{fweins}
J^{\alpha, p}(v)=\frac{\Big (\int_{\mathbb R}|D^{\frac{\alpha}{2}}v|^2\Big)^{\frac{p}{2\alpha}}{\Big (\int_{\mathbb R}|v|^2\Big)^{\frac{p}{2\alpha}(\alpha-1) + 1}}}{\int_{\mathbb R}|v|^{p+2}},
\end{equation}
then we say that $Q\in H^{\alpha/2}(\mathbb R)$ is a {\bf{ground state solution}} of equation (\ref{dground}). Here, $0< \alpha< 2$  and $0<p<p_{max}(\alpha)$, and where the critical exponent $p_{max}(\alpha)$ is defined as
\begin{equation}\label{exp}
p_{max}(\alpha)\equiv \left\{\begin{array}{lll}
\frac{2\alpha}{1-\alpha},\qquad \text{for}\;\;0<\alpha<1,\\
+\infty,\qquad \text{for}\;\; 1\leqq \alpha <2.
\end{array}\right.
\end{equation}
\end{defi}

 From Frank\&Lenzmann (Proposition 1.1 and Theorem 2.2 in \cite{FL}) there is a unique (modulo translation) ground state solution for  (\ref{dground}). For  $p<p_{max}(\alpha)$ and $\alpha\in (0,1)$, it is so-called  the case $H^{\alpha/2}$-subcritical, because of this condition on $p$ is necessary to have the existence of solutions for (\ref{dground}) (see the analysis above). Thus, via a scaling argument, we obtain that equation in (\ref{equa}) has a unique ground state solution, denoted by $Q_c$.
 Moreover, we have the following regularity and decay properties for $Q_c$: $Q_c\in H^{\alpha +1}(\mathbb R)\cap C^\infty(\mathbb R)$,  
\begin{equation}\label{decay}
\frac{C_0}{1+|x|^{\alpha +1}}\leqq |Q_c(x)|\leqq \frac{C}{1+|x|^{\alpha +1}},\;\;\;\; |xQ_c' (x)|\leqq \frac{C}{1+|x|^{\alpha +1}},\qquad \text{for all}\;\; x\in \mathbb R,
\end{equation}
with some constants $C\geqq C_0>0$ depending of $\alpha, p$ and $c$. 

The  study of stability properties for the solitary wave profile $\varphi$ in (\ref{equa}) for the case $\alpha\geqq 1$ is well developed.  Indeed, in few words, there are two useful lines of exploration for studying this relevant property in the vicinity of the wave $\varphi$.  First, we have a global variational characterization of solutions of (\ref{equa}) such that a profile $\varphi$ satisfying that $\varphi>0$ on $\mathbb R$, $\varphi$ even and $\varphi'<0$ on $(0, +\infty)$ can be seen as the infima of the constrained-mass energy minimizer
\begin{equation}\label{variational}
J=inf\{E(v): v\in H^{\alpha/2}(\mathbb R) \;\;\text{and}\;\; \int_{\mathbb R} v^2 dx=\lambda\}
\end{equation}
with $E$ being the conservation-energy functional 
\begin{equation}\label{energy}
E(v)=\frac12 \int_{\mathbb R} |D^{\alpha/2}v|^2 - \frac{2}{(p+1)(p+2)} v^{p+2}dx.
\end{equation}
We recall, since $\alpha \geqq 1$, the Sobolev embedding $H^{\alpha/2}(\mathbb R)\hookrightarrow   L^{p+2}(\mathbb R)$  ensures that the functional $E$ is well-defined for any $p\geqq 0$ and the infimum in (\ref{variational}) will satisfy $-\infty<J<0$ exactly for $p<2\alpha$  (the so-called $L^2$-subcritical case). Thus, the Concentration-Compactness Method will work  very well for obtaining both existence and stability properties of $\varphi$. More exactly, in this case we obtain the global stability property of the nonempty set of minimizer $\mathcal G$ associated to the variational problem (\ref{variational}), 
\begin{equation}\label{energy1}
\mathcal G=\{v\in H^{\alpha/2}(\mathbb R): E(v)=J\;\; \text{and}\;\; \int_{\mathbb R} v^2 dx=\lambda\}.
\end{equation}
 Thus, via a scaling argument and from the uniqueness results of the ground state solutions $Q_c$ of (\ref{equa}) for $1\leqq \alpha\leqq 2$ (see Remark 2.1 in \cite{FL}), we obtain for a specific choice of $\lambda$ in (\ref{energy1}) that
\begin{equation}\label{energy2}
\mathcal G=\{Q_c(\cdot+y): y\in \mathbb R\}\equiv \Omega_{Q_c},
\end{equation}
where $\Omega_{Q_c}$ is called the orbit generated by $Q_c$ via the basic symmetry of  translations associated to the model (\ref{0gfkdv}). For the case $p\geqq 2\alpha$, it is well known that the profile $\varphi$ is nonlinearly unstable (see Bona\&Souganidis\&Strauss \cite{BSS} for the case $p> 2\alpha$ ($\alpha\geqq 1$), Martel\&Merle \cite{MM2}-\cite{MM3} for $\alpha=2$, $p=4$, and Merle\&Pilod \cite{MP} for $\alpha=1$, $p=2$). We note that by using a variational approach, it is also possible to obtain   the instability result in the $L^2$-supercritical case $p> 2\alpha\geqq 2$  (see Angulo's book, Chapter 10, \cite{angulo4})

An similar approach of  stability for the orbit $\Omega_{Q_c}$ in the case $\frac12<\alpha<1$ and $ p<2\alpha$ ($p=1$), it has been established recently by Linares\&Pilod\&Saut in \cite{fds}. The Concentration-Compactness Method was applied successfully to the minimizer problem  in (\ref{variational}) and again the property $-\infty<J<0$ is necessary for the stability result. In this point, it is also worth noting that for $\alpha=1/2$, $J=0$. Indeed, since the ground state for (\ref{equa}) with $\alpha>1/3$ is characterized (via a scaling) as the solution of the minimization problem $J^{\frac12, 1}$ in (\ref{minimi}), we obtain from $J^{\frac12, 1}(Q_c)=\frac13 \|Q_c\|$  the sharp inequality
\begin{equation}\label{G-N}
\frac13 \int_{\mathbb R}|v|^{3} dx\leqq \frac{\|v\|}{\|Q_c\|}\int_{\mathbb R}|D^{\frac14}v|^2dx.
\end{equation}
Thus, for the restriction $\|v\|=\|Q_c\|$ we obtain immediately that $E(v)\geqq 0$ and $E(Q_c)=0$. Moreover, from (\ref{G-N}) it follows the following main property:
\begin{equation}\label{posiener}
\text{if}\;\;\; \|v\|\leqq \|Q_c\|\qquad \text{then}\qquad E(v)\geqq 0.
\end{equation}

 We recall that the later result is similar to  that for  $\alpha\geqq 1$ and $2\alpha=p$, namely, $J^{\alpha, 2\alpha}(Q_c)=0$, the so-called $L^2$-critical case (we note that, for $\alpha\in [1,2]$ and $|u|^{2\alpha}u_x$ as the nonlinear part in (\ref{0gfkdv}), recently Kenig\&Martel\&Robbiano in \cite{KMR} have proved for $\alpha$ close to 2, solutions of negative energy $E$ close to the ground state blow up in finite or infinite time in the energy space $H^{\frac{\alpha}{2}}(\mathbb R)$). The case $\alpha=\frac12$ is so-called the {\it critical case} for the fKdV model (\ref{0gfkdv}). 
 
 From the recently numerical study in Klein\&Saut in \cite{KS},  the simulations showed a possible blow-up phenomenon of the associated solutions for (\ref{0gfkdv}) with an initial data $u_0$ of negative energy ($E(u_0)<0$) and therefore with a mass larger that the ground state mass $Q_c$ ($\|Q_c\|<\|u_0\|$) (see Fig. 10 in \cite{KS}). Here we will show in Theorem \ref{t4.9} below, that in fact for this regimen of $\alpha$ we have a kind of ``{\it stability of the blow-up}'' near to the possible unstable ground state solutions and so checking one of the conjectures emerging of the numerical findings in  \cite{KS}.

The second approach for  an analysis of orbital stability is that of local type, more exactly, it is fixed a solitary wave profile $\varphi_c$ of (\ref{equa}) and we study the behavior of the  flow associated to (\ref{0gfkdv}) in a neighborhood of the orbit $\Omega_{\varphi_c}$. The main  property of the energy $E$ to be obtained in this case is the following: 
\begin{equation}\label{mainineq0}
 \left\{\begin{array}{lll}
\text{There are}\;\; \delta>0\;\;\text{and}\;\;\beta_0 >0\;\;\text{such that}\\
\\
E(v)-E(\varphi_c)\geqq \beta_0[d(v;\Omega_{\varphi_c})]^2\\
\\
\text{for}\;\;  d(v;\Omega_{\varphi_c})<\delta\;\;\text{and}\;\; F(v) = F(\varphi_c),
\end{array}\right.
\end{equation}
with $F(v)=\frac12 \int v^2dx$ and $d(v;\Omega_{\varphi_c})=\inf_{y\in \mathbb R}\|v-\varphi_c(\cdot+y)\|_{H^{\frac{\alpha}{2}}}$. So, from (\ref{mainineq0}), the  continuity of the functional $E$ and of the flow $t\to u(t)$, we obtain immediately the stability property of $\Omega_{\varphi_c}$ by initial perturbations in the manifold
\begin{equation}\label{manifold}
\mathcal M=\Big\{v: \int_{\mathbb R} v^2dx =  \int_{\mathbb R} \varphi_c^2dx\Big\}.
\end{equation}
The stability for general perturbations of $\Omega_{\varphi_c}$ can be obtained via the existence of a regular curve of solitary waves, $c\to \varphi_c$.

Now, a way for  obtaining (\ref{mainineq0}) is to use Taylor's theorem and so the analysis is reduced to study the quadratic form $\langle \mathcal L_c f, f\rangle$ on the tangent space to the manifold $\mathcal M$ at the point $\varphi_c$, $T_{\varphi_c} \mathcal M$. Here $\mathcal L_c$ represents the second variation of the action $S(v)=E(v)+c F(v)$ at the point $v=\varphi_c$, namely, the unbounded self-adjoint operator 
\begin{equation}\label{self}
S'(\varphi_c)\equiv \mathcal L_c= D^\alpha +c-\varphi_c^p
\end{equation} 
with domain $\mathcal D (\mathcal L_c)= H^\alpha(\mathbb R)$. Thus, it is well known that proving the inequality
\begin{equation}\label{liapunov}
\langle \mathcal L_c f, f\rangle\geqq \beta_1 \|f\|_{H^{\frac{\alpha}{2}}} \qquad \text{for every}\;\; f\in T_{\varphi_c} \mathcal M \cap Ker(\mathcal L_c)^\bot
\end{equation} 
for $\beta_1>0$ and $ Ker(\mathcal L_c)$ representing the kernel of $\mathcal L_c$, we obtain the key inequality (\ref{mainineq0}). The direct check of condition (\ref{liapunov}) is in general extremely inconvenient, because no requirement is directly related to the number (counting multiplicity) of negative eigenvalue of $\mathcal L_c$ (it which will be denoted henceforth by $n(\mathcal L_c)$, in other words, the Morse index of $\mathcal L_c$). Moreover, in general this operator has a nontrivial negative eigenspace. Indeed, for $\varphi_c$ being a positive solitary wave solution we obtain immediately $\langle \mathcal L_c \varphi_c , \varphi_c\rangle <0$ and so the Mini-Max principle implies $n(\mathcal L_c)\geqq 1$. The works in Benjamin \cite{Be2}, Weinstein \cite{W1}-\cite{W2} and Grillakis\&Shatah\&Strauss \cite{grillakis1} finesses this difficulty and provides a nice test that guarantees when (\ref{liapunov}) is satisfied. More exactly, we suppose that $n(\mathcal L_c)=1$, $Ker(\mathcal L_c)=[\frac{d}{dx} \varphi_c]$ and the remainder of the spectrum of $\mathcal L_c$ is positive and bounded away from zero. Then, the strictly increasing property of the mapping $c\to  \int_{\mathbb R} \varphi_c^2dx$ will imply inequality (\ref{liapunov}) and so the stability property of $\Omega_{\varphi_c}$ follows from (\ref{mainineq0}). 

Next, we call the attention about the assumption of the existence of a $C^1$-mapping $c\to \varphi_c$ of solitary waves. If we assume this condition hold for every $c>0$ and by considering the new variable 
$$
\phi(x)=c^{-\frac{1}{p}}\varphi_c(c^{-\frac{1}{\alpha}} x),
$$
we see that $\phi$ will be a solution of 
\begin{equation}\label{c=1}
D^\alpha \phi + \phi - \frac{1}{p+1} \phi^{p+1}=0.
\end{equation} 
Note the independence of $\phi$ with regard to the wave-speed $c$. Therefore,
\begin{equation}\label{1d(c)}
\frac{d}{dc} \int_{\mathbb R} \varphi_c^2dx=\|\phi\|^2  \frac{d}{dc} c^{\frac{2}{p}-\frac{1}{\alpha}}=\Big(\frac{2}{p}-\frac{1}{\alpha}\Big )c^{\frac{2}{p}-\frac{1}{\alpha}-1}\|\phi\|^2.
\end{equation}
Therefore,
\begin{equation}\label{2d(c)}
\frac{d}{dc} \int_{\mathbb R} \varphi_c^2dx>0 \Leftrightarrow p<2\alpha.
\end{equation}
Thus we see that condition in (\ref{2d(c)}) is the same  imposed for obtaining a minimum of  the variational problem (\ref{variational})  at least for $\alpha>\frac12$, and therefore it is not a technical condition for the method works !.

Next, if we consider that the curve $c\to  \varphi_c$ has a sufficiently regularity, then differentiating (\ref{equa}) with regard to the variable $c$, we obtain that 
\begin{equation}\label{3d(c)}
\mathcal L_c\Big( -\frac{d}{dc} \varphi_c \Big )=\varphi_c.
\end{equation}
Now, if for some $\psi \in D(\mathcal L_c)$ we have that $\mathcal L_c\psi=\varphi_c$, then from (\ref{3d(c)}) it follows that
$$
\mathcal L_c\Big( \frac{d}{dc} \varphi_c +\psi \Big )=0.
$$
Hence, if we suppose that $Ker(\mathcal L_c)=[\frac{d}{dx} \varphi_c]$ then $\frac{d}{dc} \varphi_c +\psi =\theta \frac{d}{dx} \varphi_c$, for $\theta\in \mathbb R$, and therefore
$$
\langle \psi, \varphi_c \rangle =-\frac 12 \frac{d}{dc} \int_{\mathbb R} \varphi^2_c dx.
$$
So, we have that the condition of strictly increasing of the mapping $c\to  \int_{\mathbb R} \varphi_c^2dx$ can be replaced by the condition:
\begin{equation}\label{new}
\text{if} \;\;\mathcal L_c\psi=\varphi_c,\quad \text{then}\quad \langle \psi, \varphi_c \rangle =\langle \mathcal L^{-1}_c \varphi_c, \varphi_c \rangle <0.
\end{equation}
Condition (\ref{new}) is useful in situations where it is not clear the existence of a family of solitary waves $\varphi_c$ depending smoothly on $c$ (see Albert \cite{Al1}). We recall that as $\mathcal L_c$ is a self-adjoint operator and $\varphi_c\in Ker(\mathcal L_c)^\bot$, the Fredholm solvability theorem guarantees always the existence of an element $\psi \in D(\mathcal L_c)$ such that $\mathcal L_c\psi=\varphi_c$.

Before establishing our first stability result, let us to define orbital stability for equation (\ref{0gfkdv}).  If $\varphi$ is a given solitary wave solution of (\ref{equa}); define for any $\eta>0$ the set $U_\eta\subset H^{\frac{\alpha}{2}}(\mathbb R)$ by 
$$
U_\eta=\{v\in H^{\frac{\alpha}{2}}: inf_{y\in \mathbb R} \|v- \varphi(\cdot+y)\|_{H^{\frac{\alpha}{2}}}<\eta\}.
$$

\begin{defi}\label{nonstab}
$\varphi$ is defined to be (orbitally) stable in $H^{\frac{\alpha}{2}}$ if
\begin{enumerate}
\item [(i)] there is a Banach space $Y\subset H^{\frac{\alpha}{2}}$ such that for all $u_0\in Y$, there is a unique solution $u$ of (\ref{0gfkdv}) in $C(\mathbb R; Y)\subset C(\mathbb R; H^{\frac{\alpha}{2}})$ with $u(x,0)=u_0$; and

\item [(ii)] for every $\epsilon>0$, there exists a $\delta>0$ such that for all $u_0\in U_\delta \cap Y$, the solution $u$ of (\ref{0gfkdv}) with $u(x,0)=u_0$ satisfies $u(t)\in U_\epsilon$ for all $t>0$.
\end{enumerate}

In the case $u\in C((-T^*, T^*); Y)\subset C((-T^*, T^*); H^{\frac{\alpha}{2}})$, where $T^*$ is the maximal time of existence of $u$, the property of stability is called conditional.
\end{defi}

Our first theorem of orbital stability for  (\ref{0gfkdv}) with a ``lower'' dispersion  (more exactly, of conditional type for $\alpha\in (\frac12,1)$, see Remark \ref{rnon} below)  is the following (\cite{Al1}, \cite{BSS}, {\cite{fds}).

\begin{teo}\label{nonlinear}$[${\bf{nonlinear stability of the ground state}}$]$ Let $\frac12 < \alpha <2$ and $0<p<p_{max}(\alpha)$. Then  the ground state solution $Q_c$ for equation (\ref{equa}) is $H^{\frac{\alpha}{2}}(\mathbb R)$-stable by the flow of equation (\ref{0gfkdv}) for $p<2\alpha$.
\end{teo}

Now, with regard to the stability (linear instability) properties of the solitary waves for (\ref{equa}) in the  case $\alpha \in (\frac13, \frac12]$ and $p=1$ in (\ref{0gfkdv}), by considering the new variable
$$
w(x,t)=u(x+ct,t)-\varphi_c(x)
$$
into the  fKdV equation and using equation (\ref{equa}) satisfied by $\varphi_c$, one finds that $w$ satisfies the nonlinear equation
\begin{equation}\label{fkdv2}
(\partial_t-c\partial_x)w+\partial_x(\varphi_c w-D^\alpha w+O(\|w\|^2))=0.
\end{equation}
As a leading approximation for small perturbation, we replace (\ref{fkdv2}) by its linearization around $\varphi_c$, and hence obtain the linear equation
\begin{equation}\label{lifkdv}
(\partial_t-c\partial_x)w+\partial_x(\varphi_cw-D^\alpha w)=0.
\end{equation}
Since $\varphi_c$ depends on $x$ and but not $t$, the equation (\ref{lifkdv}) admits treatment by separation of variables, which leads naturally to a  spectral 
problem. Seeking particular solutions of (\ref{lifkdv}) of the form $w(x,t)=e^{\lambda t} u(x)$ (so-called \textit{growing mode solution}), where $\lambda \in \mathbb C$, $u$ satisfies  the linear problem
\begin{equation}\label{spefkdv}
(\lambda- c\partial_x)u+\partial_x(\varphi_c u-D^\alpha u)=0.
\end{equation}

We can say from (\ref{spefkdv}) that the complex growth rate $\lambda$ appears as (spectral) parameter for the extended eigenvalue problem
\begin{equation}\label{specp}
\partial_x \mathcal L_cu= \lambda u,
\end{equation}
with $ \mathcal L_c$ defined in (\ref{self}) with $p=1$. If equation (\ref{specp}) have a nonzero solution $u \in D(\mathcal L_c)=H^\alpha(\mathbb R)$ then an bootstrapping  argument shows that $u\in H^s(\mathbb R)$ for all $s\geqq 1$, so that (\ref{specp}) is satisfied in classical sense.  A necessary condition for the ``stability'' of $\varphi_c$ is that there are not  points $\lambda$ with $\mbox{Re}(\lambda)>0$ (which would imply the existence  of a solution $u=u(x)$ of (\ref{specp}) that  grows exponentially in time). If we denote by $\sigma$ the ``spectrum''  of $\partial_x\mathcal{L}_c$ (namely, $\lambda \in \sigma$ if there is a $u\neq 0$ satisfying (\ref{specp})), the later discussion
suggests the utility of the following definition:

\begin{defi} (linear stability and instability) \label{defspe} A solitary wave solution $\varphi_c$ of the fKdV equation (\ref{0gfkdv}) is said to be linearly stable if $\sigma\subset i\mathbb R$. Otherwise (i.e., if $\sigma$ contains point with $\mbox{Re}(\lambda)>0$) $\varphi_c$ is  linearly unstable.
\end{defi}

We recall that as (\ref{lifkdv}) is a real Hamiltonian equation, it forces certain elementary symmetries on the spectrum of $\sigma$, more exactly, $\sigma$ will be symmetric with respect to reflection in the real and imaginary axes. Therefore, it implies that exponentially growing perturbation are always paired with exponentially decaying ones. It is the reason by which was only required in  Definition \ref{defspe} that the spectral parameter $\lambda$ satisfies that $\mbox{Re}(\lambda)>0$.

An similar spectral problem to (\ref{specp}) for traveling wave solutions (solitary or periodic)  has been the focus of many research studies in the last years, see Grillakis\&Shatah\&Strauss \cite{grillakis2}, Lopes \cite{lo}, Lin  \cite{lin}, Kapitula\&Stefanov \cite{KSt}, among others.

Our linearized instability result for the  fKdV equation (\ref{0gfkdv})
is the following: 

\begin{teo}\label{teo2}$[${\bf{Linear instability criterium for fKdV equations}}$]$
Let $c\to  \varphi_c \in H^{\alpha+1} (\mathbb R)$ be a smooth curve of positive solitary wave solutions to
equation $(\ref{equa})$ with $\alpha \in (\frac13, \frac12)$, $p=1$. The wave-speed $c$ can be considered over some nonempty interval $I$, $I\subset (0, +\infty)$. We assume that the self-adjoint operator $\mathcal L_c= D^\alpha + c-\varphi_c$ with domain $D(\mathcal L_c)= H^ \alpha(\mathbb R)$ satisfies
\begin{equation}\label{eq14}
Ker(\mathcal{L}_c)=\Big [\frac{d}{dx}\varphi_c\Big ].
\end{equation}
Denote by $n(\mathcal{L}_c)$ the number (counting multiplicity)
of negative eigenvalues of the operator $\mathcal{L}_c$.  Then there is a purely growing mode
$e^{\lambda t}u(x)$ with $\lambda>0$, $u\in H^s(\mathbb R) -\{0\}$, $s\geqq 0$, to the linearized equation $(\ref{lifkdv})$ if one
of the following two conditions is true:
\begin{enumerate}
\item[(i)] $n(\mathcal{L}_c)$ is even and $\frac{d}{dc}\langle \varphi_c, \varphi_c \rangle >0$.
\item[(ii)] $n(\mathcal{L}_c)$ is odd and $\frac{d}{dc}\langle \varphi_c, \varphi_c \rangle <0$.
\end{enumerate}

\end{teo}

The proof of the instability criterium established in Theorem \ref{teo2} is based in the compactness of some specific commutators associated to the family $\mathcal A^\lambda$ defined in (\ref{A}) below, and the  analytic and  asymptotic perturbation theory for linear operators. These approach can be also applied to the general model (\ref{gfkdv}) with a linear operator $\mathcal M$ of ``lower dispersion'' under some specific conditions about the symbol $\beta$. Also,  our approach can be extend to the  case of periodic traveling waves solutions associated to the model (\ref{gfkdv}) (a work in progress).

\begin{obs}\label{rem0}
\begin{itemize}
\item[1)] The conditions $(i)-(ii)$ in Theorem \ref{teo2} are similar to that obtained in Lin \cite{lin}, case $\alpha\geqq 1$, and  in Kapitula\&Stefanov \cite{KSt}-Corollary 16.

\item[2)] Our approach provides the existence of a nonzero solution $u\in D(\mathcal{L}_c)=H^\alpha(\mathbb R)$ satisfying the eigenvalue problem (\ref{specp}),  via a different approach than that given in Kapitula\&Stefanov \cite{KSt} and Pelinovsky \cite{P}.

\item[3)] If there is $\psi\in D(\mathcal L_c)$ such that $\mathcal L_c \psi=\varphi_c$,    the conditions $(i)-(ii)$ in Theorem \ref{teo2} can be changed by
\begin{enumerate}
\item[(i)] $n(\mathcal{L}_c)$ is even and $\langle \psi, \varphi_c \rangle <0$.
\item[(ii)] $n(\mathcal{L}_c)$ is odd and $\langle \psi, \varphi_c \rangle >0$.
\end{enumerate}

\item[4)]  The former criterium 3) is very useful when we do not have in hands a smooth curve $c\to \varphi_c$ of solitary waves. 

\end{itemize}
\end{obs}

As consequence of Theorem \ref{teo2} we obtain the following stability result for the ground state solutions of  equation (\ref{equa}) (see  \cite{KSt}).

\begin{coro}\label{coro3}$[${\bf{Linear instability of ground state for fKdV equations}}$]$
For $\alpha \in (\frac13, \frac12)$, $p=1$ and $c>0$, the ground state profiles $Q_c$   for (\ref{equa}) are spectrally unstable.
\end{coro}

Next, we consider $u(x,t)=\psi_c(x-ct)$ a solitary wave solution for the model  (\ref{gfkdv}). Then $\psi_c$ satisfies 
\begin{equation}\label{psi}
\mathcal M \psi_c+ c\psi_c -f(\psi_c)=0,
\end{equation}
and, similarly as in the case of model (\ref{0gfkdv}), we also have the linearized equation around $\psi_c$
\begin{equation}\label{lingkdv}
(\partial_t-c\partial_x)u+\partial_x(f'(\psi_c)u-\mathcal{M}u)=0.
\end{equation}
In order to obtain a growing mode solution of the form $e^{\lambda t}w(x)$,
$\mbox{Re}\lambda>0$, function $w$ must satisfy
\begin{equation}\label{growkdv}
(\lambda-c\partial_x)w+\partial_x(f'(\psi_c)w-\mathcal{M}w)=0.
\end{equation}

Then similarly as in the case of model  (\ref{0gfkdv}) we obtain the following
linearized instability result for the gfKdV equation (\ref{gfkdv}),  provided that the symbol $\beta$ defining the pseudo-differential operator $\mathcal M$ satisfies  for $\alpha\in (0,1)$ the following condition:
\begin{equation}\label{hbeta}
\text{if}\;\;\;\;\eta(\xi)\equiv \beta(\xi) -|\xi|^\alpha\quad \text{then}\quad \eta'\in L^2(\mathbb R).
\end{equation}\

\begin{teo}\label{teo3}$[${\bf{Linear instability criterium for gfKdV equation}}$]$
Let $c\in I\subset (0, +\infty)\to \psi_c\in H^{\alpha+1}$, $0<\alpha<1$, be a smooth curve of positive solitary wave solutions to equation $(\ref{psi})$. We assume condition (\ref{hbeta}) and that the self-adjoint operator $\mathcal N_c= \mathcal M+ c-f'(\psi_c)$ with domain $D(\mathcal N_c)= H^ \alpha(\mathbb R)$ satisfies
\begin{equation}\label{eq140}
Ker(\mathcal{N}_c)=[\frac{d}{dx}\psi_c].
\end{equation}
Denote by $n(\mathcal{N}_c)$ the number (counting multiplicity)
of negative eigenvalues of the operator $\mathcal{N}_c$.  Then there is a purely growing mode
$e^{\lambda t}w(x)$ with $\lambda>0$, $w\in H^s(\mathbb R) -\{0\}$, $s\geqq 0$, to the linearized equation $(\ref{lingkdv})$ if one
of the following two conditions is true:
\begin{enumerate}
\item[(i)] $n(\mathcal{N}_c)$ is even and $\frac{d}{dc}\langle \psi_c, \psi_c \rangle >0$.
\item[(ii)] $n(\mathcal{N}_c)$ is odd and $\frac{d}{dc}\langle \psi_c, \psi_c \rangle <0$.
\end{enumerate}

\end{teo}

The proof of Theorem \ref{teo3} follows the same lines of that for Theorem \ref{teo2}, but because of  the generality of the symbol associated with the operator $\mathcal M$ some points  in the analysis need to be treated carefully.

The analysis established above for the fKdV equation (\ref{0gfkdv}) can be extend to the fBBM equation (\ref{fBBM0}) for $\alpha\in (\frac13, 1)$. In section 4 below, we show the following stability properties associated to the ground state solutions $\Phi_c$ satisfying
\begin{equation}\label{sbbm1}
D^\alpha \Phi_c+ \Big(1-\frac1c\Big) \Phi_c -\frac1c \Phi^2_c=0,\qquad c>1.
\end{equation}

\begin{teo}\label{nonlinearfbbm1}$[${\bf{nonlinear stability of the ground state for the fBBM}}$]$ Let $\frac13 < \alpha <1$ and $c>1$. Then,  the ground state solution $\Phi_c$  of (\ref{sbbm1}) is $H^{\frac{\alpha}{2}}(\mathbb R)$-stable by the flow of equation (\ref{fBBM0}) provided $\alpha\in [\frac12, 1)$ and $c>1$, and for $\alpha\in (\frac13, \frac12)$ and $c>c_0$. Here $c_0$ is given by 
$$
c_0=\frac{2+\sqrt{2(3\alpha-1)}}{6\alpha}.
$$
\end{teo}

Our  nonlinear stability results for the fBBM equation extends and complements those in Linares\&Pilod\&Saut \cite{fds}, in the sense that we show stability of the orbit $\Omega_{\Phi_c}=\{\Phi_c(\cdot+y): y\in \mathbb R\}$ for $\alpha=\frac12$ and $c>1$, and  for $\frac13<\alpha<\frac12$ with the specific restriction
on the wave velocity $c$. It also confirms the numerical simulations in Klein\&Saut \cite{KS} about the stability of the solitary waves in this regimen of $\alpha$'s. Similarly as in the case of the fKdV equation, the statement of orbital stability in Theorem \ref{nonlinearfbbm1} is a conditional one (see Remark \ref{refbbm} below).

\begin{teo}\label{linsfbbm1}$[${\bf{Linear instability of ground state for fBBM equations}}$]$
For $\alpha \in (\frac13, \frac12)$ and $c\in (1, c_0)$, the ground state profiles $\Phi_c$   for (\ref{sbbm1}) are linearly unstable.
\end{teo}

\indent Our paper is organized as follows. In Section 2 we present the proof of the  linear instability criterium in Theorem \ref{teo2} for the fKdV model (\ref{0gfkdv}) and of the criterium for  the general dispersive equation (\ref{gfkdv}) in Theorem \ref{teo3}. In Section 3, we prove our ``stability of the blow-up'' for the critical fKdV equation (\ref{0gfkdv}) ($\alpha=\frac12$, $p=1$). In the final section 4, we prove the results of nonlinear stability and linear instability for the fBBM equation (\ref{fBBM0}) (Theorems \ref{nonlinearfbbm1} and \ref{linsfbbm1}).
\vskip0.2in

\noindent \textbf{Notation.} We will denote $| \cdot\ |_p$ the norm in the Lebesgue space $L^p(\mathbb R)$, $1 \leqq
p \leqq \infty$ and $\|\cdot\|_s$ the norm in the Sobolev space $H^s(\mathbb R)$, $s \in \mathbb R$. For $X, Y$ Banach spaces, $B(X;Y)$ represents the set of bounded linear operators $F:X\to Y$. $[A, B]=AB-BA$ represents the commutator of the operators $A$ and $B$. $\rho(A)$ will represent the resolvent of the linear operator $A$.

\section{Nonlinear stability and Linear instability for the fKdV equation}

This section is devoted to show Theorem \ref{nonlinear} and Theorem \ref{teo2} established in the introduction. The proof of the nonlinear stability of the ground state is a consequence immediate of Grillakis {\it et.al} \cite{grillakis1} and Frank\&Lenzmann \cite{FL} results. For the linear instability result we extend the theory of Vock\&Hunziker in \cite{VH} about the stability of Schr\"odinger eigenvalue problems to the  study of linear instability of solitary waves solutions for the fKdV type model in (\ref{0gfkdv}) with a ``lower dispersion''.  In particular, we recover the linear instability results in \cite{KSt} for the ground state solutions associated to the equation (\ref{equa})  with $p=1$ and $\alpha \in (\frac 13, \frac12)$. Our analysis is also extended to the general lower-dispersion models (\ref{gfkdv}).

\subsection{Nonlinear stability of ground state for the fKdV equation}

In the following we show Theorem \ref{nonlinear}, we recall that in the literature it result of stability has been showed by difference methods (\cite{Al1}, \cite{BSS}, \cite{fds}). Here, for completeness of the exposition, we show this in a unified way for $\alpha\in (\frac12, 2)$.

\begin{proof}$[${\bf{Proof of Theorem \ref{nonlinear}}}$]$  Let $Q$ be the ground state solution for (\ref{dground}), namely, $ Q=Q(|x|)>0$, satisfies
\begin{equation}
D^\alpha Q+Q-Q^2=0,
\end{equation}
and a  minimum for the functional $J^{\alpha, 1}$ in (\ref{minimi}). Thus we obtain that the  self-adjoint operator,
\begin{equation}
\mathcal L_1= D^\alpha +1-2Q
\end{equation}
satisfies the so-called nondegeneracy property, namely,  $Ker(\mathcal L_1)=[\frac{d}{dx}Q]$. Moreover, since $
\langle \mathcal L_1 Q, Q\rangle\leqq 0$ and for $\eta\in C^{\infty}_0(\mathbb R)$
$$
\langle \mathcal L_1 \eta, \eta\rangle\geqq 0,\qquad \text{for all}\;\; \eta\bot Q^2, 
$$
we have  $n(\mathcal L_1)=1$ (see \cite{FL}). Now, for $R\equiv\alpha Q+x Q' \in L^2 (\mathbb R)$ (see (\ref{decay})) follows $\mathcal L_1 R=-\alpha Q$ (at least in the distributional sense). Thus a bootstrapping argument shows that $R\in H^{\alpha+ 1}(\mathbb R)$ and so $R\in D(\mathcal L_1)=H^{\alpha}(\mathbb R)$.

Next,  for any real number $\theta \neq 0$, define the dilation operator $T_\theta$ by $(T_\theta f)(x)= f(\theta x)$. Then,  via the elementary scaling $Q_c(x)=2c Q(c^{1/\alpha} x)$ and the relation $D^\alpha (T_\theta f)(x)= \theta^\alpha D^\alpha f (\theta x)$, we can show that for $\theta= c^{1/\alpha}$ we obtain that $Q_c$ satisfies
\begin{equation}
D^\alpha Q_c +c Q_c-\frac12 Q_c^2=0,
\end{equation}
and so, we obtain its linearized operator 
\begin{equation}
\mathcal L_c= D^\alpha +c-Q_c.
\end{equation}
Now, the relation $\mathcal L_c= c T_\theta \mathcal L_1 T^{-1}_\theta$ implies with $\theta =c^{1/\alpha}$ that $spec(\mathcal L_c)=\{c r: r\in spec(\mathcal L_1)\}$ and therefore $\mathcal L_c$ and $\mathcal L_1$ have the ``{\it same structure}''. Thus, $\psi$ is an eigenfunction of $\mathcal L_1$ with eigenvalue $\lambda$ if and only if $T_\theta\psi$ is an eigenfunction of $\mathcal L_c$ with eigenvalue $c\lambda$. Then, we conclude immediately that $n(\mathcal L_c)=1$ and  $Ker(\mathcal L_c)=[\frac{d}{dx}Q_c]$. Thus, since $R_c=\alpha Q_c+ x Q'_c\in D(\mathcal L_c)$ with $\mathcal L_c R_c=-\alpha cQ_c$ (where we are used $D^\alpha(x Q')=\alpha D^ \alpha Q+ xD^\alpha Q'$) we obtain
\begin{equation}
 \langle \mathcal L^{-1}_cQ_c, Q_c\rangle= -\frac{1}{\alpha c} \langle R_c, Q_c\rangle =\|Q_c\|^2 \Big [\frac{1}{2\alpha c}-\frac{1}{c}\Big] <0,
\end{equation}
where we use integration by parts ($xQ_c^2(x)\to 0$ as $|x|\to +\infty$) and $\alpha>\frac12$. Hence, from regularity properties of the curve $c\to Q_c$ (see proof of Corollary \ref{coro3} below)  
and from the Lyapunov property of the energy $E$ in (\ref{mainineq0}) we finish the proof.
\end{proof}

\begin{obs}\label{rnon} The statement in Theorem \ref{nonlinear} deserves to be clarified at  least in some points with regard to  the Cauchy problem.
\begin{enumerate}
\item [(1)] For $2>\alpha> 1$ the solutions of the Cauchy problem are global in $H^{\frac{\alpha}{2}}(\mathbb R)$ and so the stability result is not conditional (see Definition \ref{nonstab}).  Indeed, by using a similar strategy to that in the proof of Theorem 1 in Kenig\&Martel\&Robbiano \cite{KMR} we obtain local well-posedness of the model in (\ref{0gfkdv}) for every initial data $u_0\in H^{\frac{\alpha}{2}}(\mathbb R)$  for $p<2\alpha$  (in \cite{KMR} was studied the case of the critical-nonlinearity $|u|^{2\alpha} u_x$). Moreover, the conservation of the energy $E$ in (\ref{energy}) and the charge $F(u)=\int_{\mathbb R} u^2dx$ by the flow of  (\ref{0gfkdv}), together with an application of the Gagliardo-Nirenberg type inequality (see (\ref{fweins}))
\begin{equation}\label{GN}
\int_{\mathbb R}|u|^{p+2}\leqq C_{\alpha, p} \Big(\int_{\mathbb R}|D^{\frac{\alpha}{2}}u|^2\Big)^{\frac{p}{2\alpha}}\Big(\int_{\mathbb R}|u|^2\Big)^{\frac{p}{2\alpha}(\alpha-1) + 1},
\end{equation}
it gives us exactly for $p<2\alpha$ the ``{\it a priori}'' estimative, 
\begin{equation}\label{apriori}
\begin{array}{lll}
\|D^{\frac{\alpha}{2}}u\|^2&\leqq  E(u_0)+ C_{\alpha, p}\|D^{\frac{\alpha}{2}}u\|^{\frac{p}{\alpha}}\|u\|^\beta\\
&\leqq C(\|u_0\|_{\frac{\alpha}{2}})+ \frac{D_{\alpha, p}}{r} \|u_0\|^{pr+2} + \frac{p}{2\alpha}\|D^{\frac{\alpha}{2}}u\|^2
\end{array}
\end{equation}
with $\beta= \frac{p}{\alpha}(\alpha-1) +2$ and $r=\frac{2\alpha}{2\alpha-p}$. Thus, global-posedness of the initial valued problem for (\ref{0gfkdv})  follows in the  energy space $H^{\frac{\alpha}{2}}(\mathbb R)$.

\item [(2)] In Herr\&Ionescu\&Kenig\&Koch \cite{HIKK} was showed that for $\alpha\in (1,2)$ the solutions of the Cauchy problem for (\ref{0gfkdv}) with $p=1$ are global well-posed in $L^2(\mathbb R)$.

\item [(3)] The case $\alpha=1$ $($ $p=1$, that is, the Benjamin-Ono equation$)$ was showed   to be  globally well-posed in $H^{s}(\mathbb R)$ for $s\geqq 0$ by Ionescu\&Kenig \cite{IK}.

\item [(4)] The case $\alpha\in (\frac12, 1)$ is more delicate with regard to the local and global well-posedness problem. In Saut \cite{S} was proved that (\ref{0gfkdv})  admit global weak solutions (without uniqueness) in the space $L^\infty(\mathbb R; H^{\frac{\alpha}{2}}(\mathbb R))$ and global weak solutions in $L^\infty(\mathbb R; L^2(\mathbb R)) \cap L^2_{loc}(\mathbb R; H^{\frac{\alpha}{2}}_{loc}(\mathbb R))$ in Ginibre\&Velo \cite{GV1}-\cite{GV2}. 

\item [(5)] The best known result of local well-posedness for (\ref{0gfkdv}) has been established by Linares\&Pilod\&Saut \cite{fds} in $ H^s(\mathbb R)$, for $s>s_\alpha\equiv \frac32-\frac{3\alpha}{8}>\frac{\alpha}{2}$, for $\alpha \in (0,1)$. It which does not allow to globalize the solution using conservation laws.

\item [(6)] The problem to prove local well-posedness in $H^{\frac{\alpha}{2}}(\mathbb R)$ in the case 
$\alpha\in [\frac12, 1)$, which would imply global well-posedness by using the conserved quantities $E$ and $F$, is still open.

\item [(7)]  Therefore, the statement of stability in Theorem \ref{nonlinear} for $\alpha\in (\frac12, 1)$ is one of {\bf{conditional}} type by the Definition \ref{nonstab}, where we  have used  $Y=H^s(\mathbb R)$, $s>s_\alpha>\frac12$, and so for all $\epsilon>0$ there is a $\delta>0$ such that if $u_0\in  H^s(\mathbb R) \cap U_\delta $,  then $u(t)\in U_\epsilon$, for all $t\in (-T_s, T_s)$, where $T_s$ is the maximal time of existence of $u$ satisfying $u(0)=u_0$.

\end{enumerate}
\end{obs}

\subsection{linear instability criterium for the fKdV equation}

In order to illustrate the strategy for obtaining a \textit{growing mode solution} of (\ref{lifkdv}) with the form $w(x,t)=e^{\lambda t}u(x)$ and $\mbox{Re}\lambda>0$, we can see the eigenvalue problem (\ref{spefkdv}) for $\lambda$ and $u$ rewrite in the form
\begin{equation}\label{rewrite}
cu + \frac{c\partial_x}{\lambda-c\partial_x}(\varphi_c u- D^\alpha u)=0.
\end{equation}
Here the expression $\frac{\partial_x}{\lambda-c\partial_x}$ is a
notation for the well-defined  linear operator
$\partial_x(\lambda-c\partial_x)^{-1}$ and $\varphi_c$ is {\it any positive solitary wave solutions for 
(\ref{equa})}. Thus, if we consider the following family of closed linear operators
$\mathcal{A}^{\lambda}:H^{\alpha}(\mathbb R)\longrightarrow L^2(\mathbb R)$, with $ \mbox{Re} \lambda
>0$, 
\begin{equation}\label{A}
\mathcal{A}^{\lambda}v\equiv cv + \frac{c\partial_x}{\lambda-c\partial_x}(\varphi_c v- D^\alpha v), 
\end{equation}
it follows immediately that the solution of the eigenvalue problem (\ref{rewrite}) is reduced to find $\lambda\in\mathbb{C}$ with $\mbox{Re}\lambda>0$ such that the
operator $\mathcal{A}^{\lambda}$ possesses a {\it nontrivial kernel}.  Now,  we note that from the analyticity of the resolvent
associated to the operator $\partial_x$, $\lambda \in \mathcal S \to
(\lambda-c\partial_x)^{-1}$, for $\mathcal S=\{z\in \mathbb C:\mbox{Re} z >0\}$,
 we obtain that the mapping $\lambda\in \mathcal S \to \mathcal{A}^{\lambda}$ represents an
analytical family of operators of type-$A$ (see Kato \cite{kato1}), namely,
\begin{enumerate}

\item [1)] $D(\mathcal{A}^{\lambda})=H^{\alpha}(\mathbb R) $ for all $\lambda\in \mathcal S$,

\item[2)] for $u\in H^{\alpha}(\mathbb R) $, $\lambda\in \mathcal S
\to \mathcal{A}^{\lambda}u$ is analytic in the topology of
$L^2(\mathbb R) $.
\end{enumerate}
Therefore, from classical analytic perturbation theory, all discrete eigenvalues of $\mathcal{A}^{\lambda}$ ($\mbox{Re} (\lambda)>0$) will be  stable :
{\it for $\eta$ in the discrete spectrum of $\mathcal A^{\lambda}$, there is $\delta>0$ such that for $\lambda_0\in B(\lambda;\delta)$, $\mathcal A^{\lambda_0}$ has $\eta_i(\lambda_0)$ eigenvalues close to $\eta$ with total algebraic multiplicity equal to that of $\eta$}.

 In our approach, we will find a growing mode solution for $\lambda>0$. Indeed, since we  have that
$$
\mathcal{A}^{\lambda}\longrightarrow \mathcal{L}_c\equiv D^\alpha +c - \varphi_c\;\;\;\;\;\text{as}\;\; \lambda\rightarrow0^{+},
$$
strongly in $L^2(\mathbb R)$ (see Proposition \ref{strong} below), we will use asymptotic perturbation arguments as  Vock\&Hunziker in \cite{VH}  and Lin in \cite{lin} (see also Hislop\&Sigal in \cite{hislop}) for obtaining our criterium established in Theorem \ref{teo2}.  In our analysis, it will be decisive to count the number of eigenvalues of
$\mathcal{A}^{\lambda}$ (for $\lambda$ small) in the left-half plane (for  $\lambda$ large, there is not growing modes, see Lemma \ref{lema3} below), and so we will
need to know how the zero eigenvalue of $\mathcal{L}_c$ will be
perturbed. To this end, we obtain a moving kernel
formula  (see Lemma \ref{lema6} below) which will decide whether zero jumps to the left or to the right. Thus we get the conditions $(i)-(ii)$ in Theorem \ref{teo2}.

As the structure of the proof of Theorem \ref{teo2} follows some ideas used by Lin in \cite{lin} for the case $\alpha\geqq 1$, we will only indicate the new basic differences due to  the structure of the operator $\mathcal A^\lambda$ defined in (\ref{A}) for $\alpha\in (0,1)$.  

\subsubsection{Stability of the discrete spectrum of $\mathcal L_c$, with $\alpha \in (0,1)$}

In this subsection we study the behavior of the family $\mathcal A^\lambda$ by depending of $\lambda$. In particular, we show that every discrete eigenvalue of the limiting operator $\mathcal L_c=D^\alpha + c-\varphi_c$   is stable with respect to the family $\mathcal A^\lambda$ for small positive $\lambda$. Our first result is about the strong convergence of $\mathcal A^\lambda$.

\begin{prop}\label{strong} For $\lambda>0$, the operator $\mathcal A^\lambda$ converges to $\mathcal L_c$ strongly in $L^2(\mathbb R)$ when $\lambda\to 0^+$.
\end{prop}

\begin{proof} For $\lambda>0$ and $v\in D(\mathcal A^\lambda)= D(\mathcal L_c)$ we have the relation
$$
(\mathcal A^\lambda- \mathcal L_c)v=\frac{\lambda}{\lambda-c\partial_x} (\varphi_c-D^\alpha) v.
$$
Thus, by Plancherel and  the dominated convergence theorem follows
$$
\|(\mathcal A^\lambda- \mathcal L_c)v\|^2=\int_\mathbb R \frac{\lambda^2}{\lambda^2+c^2\xi^2} |\widehat{\varphi_c}(\xi)-\widehat{D^\alpha}(\xi)|^2d\xi \to 0
$$
when $\lambda\to 0^+$.

\end{proof} 

Next, we   localized  the essential spectrum of $\mathcal A^\lambda$, $\sigma_{ess}(\mathcal A^\lambda)$. We will see that it set is situated  in the right-hand side of the complex-plane and away from the imaginary axis.
We star with the following two basic definition related to the $\sigma_{ess}(\mathcal A^\lambda)$ (see Hislop\&Sigal \cite{hislop}).

\begin{defi}\label{Zhislin} $[$Zhislin spectrum  $Z(\mathcal A^\lambda)$$]$
A Zhislin sequence for $\mathcal A^\lambda$ and $z\in \mathbb R$ is a sequence
$$
\{u_n\}\subset H^\alpha(\mathbb R),\;\;\;\|u_n\|=1,\;\;\;\; \text{supp}\;u_n\subset \{x: |x|\geqq n\}
$$
and $\|(\mathcal A^\lambda-z)u_n\|\to 0$ as $n\to +\infty$.

The set of all $z$ such that a Zhislin sequence exists for $\mathcal A^\lambda$ and $z$ is denoted by $Z(\mathcal A^\lambda)$.
\end{defi}

\begin{obs} Every  Zhislin sequence $\{u_n\}$ necessarily converges weakly to zero in $L^2(\mathbb R)$.
\end{obs}

\begin{defi}\label{Weyl} $[$Weyl spectrum  $W(\mathcal A^\lambda)$$]$
A Weyl sequence for $\mathcal A^\lambda$ and $z\in \mathbb R$ is a sequence
$$
\{u_n\}\subset H^\alpha(\mathbb R),\;\;\;\|u_n\|=1,\;\;\;\; u_n\to 0\;\;\text{weakly in}\;\; L^2(\mathbb R),
$$
and $\|(\mathcal A^\lambda-z)u_n\|\to 0$ as $n\to +\infty$.

The set of all $z$ such that a Weyl  sequence exists for $\mathcal A^\lambda$ and $z$ is denoted by $W(\mathcal A^\lambda)$.
\end{defi}

From the last two definitions we have the following result (see  \cite{hislop}).

\begin{prop}\label{prop0}
$Z(\mathcal A^\lambda)\subset W(\mathcal A^\lambda)$, $W(\mathcal A^\lambda)\subset \sigma_{ess}(\mathcal A^\lambda)$ and $\partial ( \sigma_{ess}(\mathcal A^\lambda))\subset W(\mathcal A^\lambda)$.
\end{prop}

Our main result about the $\sigma_{ess}(\mathcal A^\lambda)$ is the following one.

\begin{prop}\label{prop1}
For any $\lambda>0$, we have
\begin{equation}
\sigma_{ess}(\mathcal A^\lambda)\subset \big \{ z: Re z\geqq \frac12 c\Big\}
\end{equation}
\end{prop}

The idea of the proof of Proposition \ref{prop1} will be to see $W(\mathcal A^\lambda)=Z(\mathcal A^\lambda)$ and it will be based on the following two  lemmas.

\begin{lema}\label{lema1}
For any $\lambda>0$, we have
\begin{equation}
Z(\mathcal A^\lambda)\subset \big \{ z: Re z\geqq \frac12 c\Big\}.
\end{equation}
\end{lema}

\begin{proof} Let $z\in Z(\mathcal A^\lambda)$ and suppose $Re z<\frac12 c$. It is immediate from Fourier transform that for any  $u\in H^\alpha(\mathbb R)$ we have
$$
I_0(u)\equiv Re \langle  - \frac{c\partial_x}{\lambda-c\partial_x} D^\alpha u, u\rangle\geqq 0.
$$
Then, for any sequence $\{u_n\}\subset H^\alpha(\mathbb R)$, $\|u_n\|=1$, and satisfying $\text{supp}\;u_n\subset \{x: |x|\geqq n\}$, we have from the following trivial estimative for any $c$ and $\lambda$
$$
\Big\|\frac{c\partial_x}{\lambda+c\partial_x} u_n \Big\|\leqq \|u_n\|=1,
$$
that for $n$ large,
\begin{equation*}
\begin{array}{lll}
Re \langle (\mathcal A^\lambda-z) u_n, u_n\rangle
&\geqq I_0(u_n)+c- Re z+ Re \langle \frac{c\partial_x}{\lambda-c\partial_x}(u_n\varphi_c), u_n\rangle\\
\\
&\geqq c- Re z -Re \langle u_n\varphi_c, \frac{c\partial_x}{\lambda+c\partial_x} u_n\rangle\\
\\
& \geqq c- Re z -\text{sup}_{|x|\geqq n} |\varphi_c(x)|\geqq c- \frac12 c-\frac14 c=\frac14 c.
\end{array}
\end{equation*}
Then, since $|Re \langle (\mathcal A^\lambda-z) u_n, u_n \rangle |\leqq \|(\mathcal A^\lambda-z) u_n\|$, for all $n$, and  $(\mathcal A^\lambda-z) u_n\to 0$ as $n\to +\infty$, we obtain a contradiction.
\end{proof}

The next lemma extended  Lemma 2.3 in \cite{lin} to the case $\alpha\in (0,1)$.

\begin{lema}\label{lema2}
Given $\lambda>0$. Let $\zeta\in C_0^\infty(\mathbb R)$ be a cut-off function such that $\zeta|_{\{|x|\leqq R_0\}} =1$, for some $R_0$. Define $\zeta_d(x)=\zeta(x/d)$, $d>0$. Then, for each $d$, the operator $\zeta_d(\mathcal A^\lambda-z)^{-1}$ is compact for some $z\in \rho(\mathcal A^\lambda)$, and  there exists $C(d)\to 0$ as $d\to \infty$ such that for any $u\in C_0^\infty(\mathbb R)$,
\begin{equation}\label{conmuta0}
\|[\mathcal A^\lambda, \zeta_d]u\|\leqq C(d) (\|\mathcal A^\lambda u\|+\|u\|).
\end{equation}
\end{lema}

\begin{proof} Initially we prove that for $k>0$ sufficiently large, $-k\in \rho(\mathcal A^\lambda)$. Indeed, for $\lambda>0$ we write $\mathcal A^\lambda=D^\alpha +c +\mathcal K^\lambda$, where
\begin{equation}\label{K}
\mathcal K^\lambda=\frac{c\partial_x}{\lambda-c\partial_x}\varphi_c-
\frac{\lambda}{\lambda-c\partial_x}D^\alpha: L^2(\mathbb R)\to L^2(\mathbb R)
\end{equation}
is a bounded operator, because of the symbol of $\partial_x (\lambda-c\partial_x)^{-1}$ and $ (\lambda-c\partial_x)^{-1}D^\alpha$ are bounded (here we use that $\alpha<1$). Thus, since $A=D^\alpha +c$ is a nonnegative self-adjoint operator  and $\mathcal K^\lambda$ is a $A$-bounded operator with relative $A$-bound  equal to zero, we have that $-k\in \rho(A)$ for all $k>0$ and
$$
\|\mathcal K^\lambda (A+k)^{-1}u\|\leqq C_{\lambda, c}\|(A+k)^{-1}u\|\leqq C_{\lambda, c}\frac{1}{|k|}\|u\|.
$$
Therefore, the relation $\mathcal A^\lambda+k=[1+\mathcal K^\lambda (A+k)^{-1}](A+k)$ implies that for $k$ large, $z=-k\in \rho(\mathcal A^\lambda)$. Now, the compactness of $\zeta_d(\mathcal A^\lambda-z)^{-1}$ follows from the relation
$$
\|D^\alpha (\mathcal A^\lambda+k)^{-1} f_n\|\leqq \|[1+\mathcal K^\lambda (A+k)^{-1}]^{-1}f_n\|\leqq M\|f_n\|
$$
for all $L^2(\mathbb R)$-bounded sequence $\{f_n\}$, the local compactness of $H^\alpha (\mathbb R) \hookrightarrow L^2 (\mathbb R)$ and a Cantor diagonalization argument which imply that $(\mathcal A^\lambda+k)^{-1} f_n\to f$ in $L^2 _{loc}(\mathbb R)$ and so the sequence $\zeta_d(\mathcal A^\lambda-z)^{-1}f_n$ is convergent.

 To show the commutator estimative in (\ref{conmuta0}), we note initially that the graph norm of $\mathcal A^\lambda$ appearing at the right-side hand of (\ref{conmuta0}) is equivalent to the $\|\cdot\|_{H^\alpha}$-norm (it follows immediately from the relations $\mathcal A^\lambda=D^\alpha +c +\mathcal K^\lambda$ and $(A+k)= [1+\mathcal K^\lambda ( A+k)^{-1}]^{-1}(\mathcal A^\lambda+k)$). Now, it is not difficult to see that $\mathcal E^\lambda\equiv\frac{\lambda}{\lambda-c\partial_x}\in B(L^2(\mathbb R))$ and we have the equality
\begin{equation}\label{conmu2}
[\mathcal A^\lambda, \zeta_d]=(1-\mathcal E^\lambda)[D^\alpha, \zeta_d] + [\mathcal E^\lambda, \zeta_d](\varphi_c-D^\alpha).
\end{equation}
Next, we estimative every term at the right-hand side of (\ref{conmu2}). First, from the relation 
\begin{equation*}
\begin{split}
[\mathcal E^\lambda, \zeta_d] (\varphi_c-D^\alpha)&=\frac{1}{\lambda}\mathcal E^\lambda[c\partial_x, \zeta_d]\frac{1}{\lambda-c\partial_x}(\varphi_c-D^\alpha)\\
&=\frac{c}{\lambda d}\mathcal E^\lambda (\zeta'(x/d) \mathcal E^\lambda(\varphi_c-D^\alpha)),
\end{split}
\end{equation*}
we obtain

\begin{equation}\label{conmu3}
\begin{split}
\|[\mathcal E^\lambda, \zeta_d] (\varphi_c-D^\alpha)u\|&\leqq \frac{c}{\lambda d}
\|\zeta'\|_\infty\|\varphi_cu-D^\alpha u\|\\
&\leqq \frac{C_0}{\lambda d} (\|u\| + \|D^\alpha u\|).
\end{split}
\end{equation}
Now, since $1-\mathcal E^\lambda\in B(L^2(\mathbb R))$, we obtain from Theorem 3.3 in  Murray \cite{Murray} the estimative for $\alpha \in (0,1)$
\begin{equation}\label{conmu4}
\|(1-\mathcal E^\lambda)[D^\alpha, \zeta_d] u\|\leqq \|[D^\alpha, \zeta_d] u\|\leqq K(\alpha)\|D^\alpha \zeta_d\|_\ast \|u\|\leqq \frac{C_1}{d^\alpha}\|D^\alpha \zeta\|_\infty \|u\|,
\end{equation}
where $\|\cdot\|_\ast $ is the BMO norm and we are using the identity $(D^\alpha \zeta_d)(x)= \frac{1}{d^\alpha} (D^\alpha \zeta)(x/d)$ and the classical embedding $L^\infty(\mathbb R) \hookrightarrow BMO$ ($\|f\|_\ast \leqq 2\|f\|_\infty$). Therefore, from (\ref{conmu2})-(\ref{conmu3})-(\ref{conmu4})  we obtain the right-hand side in (\ref{conmuta0}). It finishes the proof of the Lemma.
\end{proof}

\begin{proof} $[${\bf{Proposition \ref{prop1}}}$]$. By using Theorem 10.12 in \cite{hislop} and Lemma \ref{lema2} above, we have for any $\lambda>0$   that $W(\mathcal A^\lambda)=Z(\mathcal A^\lambda)$.  Therefore, Lemma \ref{lema1} and Proposition \ref{prop0} imply Proposition \ref{prop1}. Indeed, suppose that that $z\in \sigma_{ess}(\mathcal A^\lambda)$ and $Re z<\frac12 c$. Then $z\in \mathbb C-W(\mathcal A^\lambda)\subset  \mathbb C-\partial(\sigma_{ess}(\mathcal A^\lambda))$. Therefore, $z\in Int(\sigma_{ess}(\mathcal A^\lambda))$. So, if we consider $C_z$ being the maximal non-empty open connected component of  $Int(\sigma_{ess}(\mathcal A^\lambda))$ containing the point $z$, we see that $\partial C_z \cap \{z: Re z<\frac12 c\}\neq \emptyset $. Therefore, since $\partial C_z\subset \partial(\sigma_{ess}(\mathcal A^\lambda))$ we obtain $ \partial(\sigma_{ess}(\mathcal A^\lambda)) \cap \{z: Re z<\frac12 c\}\neq \emptyset $ and so $ W(\mathcal A^\lambda) \cap \{z: Re z<\frac12 c\}\neq \emptyset $, it which is a contradiction. This finishes the proof.

\end{proof}

Next, we study the behavior of $\mathcal A^\lambda$ near infinity. We will show the non-existence of growing modes at the left-hand  side of the complex-plane for large $\lambda$ (so, since the eigenvalues of $\mathcal A^\lambda$ appear in conjugate pairs, there are not  growing modes in all for large $\lambda$).

\begin{lema}\label{lema3} There exists $\Lambda>0$, such that for $\lambda>\Lambda$, $\mathcal A^\lambda$ has no eigenvalues in $\{z: Re z\leqq 0\}$.
\end{lema}
The proof of Lemma \ref{lema3} is the same as that of Lemma 4.1 in \cite{lin} which works still for $\alpha\in (0,1)$. Next, we study the behavior of $\mathcal A^\lambda$ for small positive $\lambda$ and it is the more delicate part in the theory. The next two result are the heart of the Vock\&Hunziker theory for obtaining that every discrete eigenvalue of the limiting operator $\mathcal L_c=D^\alpha + c-\varphi_c$   is stable with respect to the family $\mathcal A^\lambda$ for small positive $\lambda$ (see Chapter 19 in Hislop\&Sigal \cite{hislop}). The following result extends those of Lin in \cite{lin} for the case $\alpha\in (0,1)$.

\begin{lema}\label{lema4} Given $F\in C_0^\infty (\mathbb R)$. Consider any sequence $\lambda_n\to 0^+$ and $\{u_n\} \subset  H^\alpha(\mathbb R)$ satisfying
\begin{equation}\label{conmuta40}
\|\mathcal A^{\lambda_n} u_n\|+\|u_n\|\leqq M_1<\infty
\end{equation}
for some constant $M_1$. Then if $w-\lim_{n\to \infty} u_n=0$, we have
\begin{equation}\label{conmuta50}
\lim_{n\to \infty}\|Fu_n\|=0
\end{equation}
and
\begin{equation}\label{conmuta60}
\lim_{n\to \infty} \|[\mathcal A^{\lambda_n}, F] u_n\|=0.
\end{equation}
\end{lema}

\begin{proof} The convergence in (\ref{conmuta50}) is immediately. Next, from the relation $\mathcal A^{\lambda_n}=D^\alpha +c+\mathcal K^{\lambda_n}$, with $\mathcal K^{\lambda_n}$ defined in (\ref{K}), we have
\begin{equation}\label{conmuta7}
[\mathcal A^{\lambda_n}, F] u_n=[D^\alpha, F] u_n +[\mathcal K^{\lambda_n}, F]u_n.
\end{equation}
Now, we show that every commutator at the right-hand side of (\ref{conmuta7}) goes to zero for $n\to \infty$.

\begin{enumerate}
\item[1)] $[D^\alpha, F] u_n\to 0$ as $n\to \infty$: for  $\mathcal C_\alpha\equiv [D^\alpha, F] $ we have that $\mathcal C_\alpha: L^2(\mathbb R)\to L^2(\mathbb R)$ is a compact operator for every $\alpha\in (0,1)$. Indeed, from Theorem 3.3 in \cite{Murray} we know that  $\mathcal C_\alpha$ is a bounded operator on $L^2(\mathbb R)$. With regard to the compactness property, we have from relation (3.18) in \cite{Murray} the representation formula
\begin{equation}\label{represe}
\mathcal C_\alpha=ik(\alpha)\int_0^\infty [F, P_t] t^{-\alpha} \frac{dt}{t},\quad\;\; \alpha\in (0,1),
\end{equation}
where $P_t=(1-t^2\partial_x^2)^{-1}$ and $k(\alpha)=(2i/\pi)sin(\pi\alpha/2)$. Now, since the symbol associated to the operator  $P_t$, $p_t(\xi)=\frac{1}{1+t^2 \xi^2}$, satisfies that $\frac{d}{d\xi}p_t(\xi)\to 0$ when $|\xi|\to \infty$, by Theorem C in \cite{cords}  we obtain that the commutator $[F, P_t]: L^2(\mathbb R)\to L^2(\mathbb R)$ is compact. Therefore $\mathcal C_\alpha$ is compact  because is the limit of a sequence of compact operators. It finishes this item.

\item[2)] $[\mathcal K^\lambda, F]u_n\to 0$ as $n\to \infty$:  We star for studying the commutator
\begin{equation}\label{conmuta8}
[\mathcal E^\lambda D^\alpha, F] u_n=  \mathcal E^\lambda D^\alpha(Fu_n)-F\mathcal E^\lambda D^\alpha u_n.
\end{equation} 
with $ \mathcal E^\lambda=\frac{\lambda}{\lambda-c\partial_x}$. Indeed, since $Fu_n\to 0$ on $L^2(\mathbb R)$ by (\ref{conmuta50}) and $ \mathcal E^\lambda D^\alpha: L^2(\mathbb R)\to L^2(\mathbb R)$ is a bounded operator for $\alpha\in (0,1)$ we have immediately that the first term of the right-hand side of (\ref{conmuta8}) goes to zero. Now, we consider $\delta>0$ such that  $\delta+\alpha<1$. Then it is not difficult to see that $\{\mathcal E^\lambda D^\alpha u_n\}$ is a bounded sequence in $H^\delta(\mathbb R)$. Therefore, by the local compact embedding of $H^\delta(\mathbb R)$ in $L^2(\mathbb R)$ and a Cantor diagonalization argument we have $\mathcal E^\lambda D^\alpha u_n\to 0$ in $L^2_{loc}(\mathbb R)$. Thus, $F\mathcal E^\lambda D^\alpha u_n\to 0$ in $L^2(\mathbb R)$. 

Now, we study the commutator
\begin{equation}\label{conmuta9}
[ \frac{c\partial_x} {\lambda-c\partial_x}\varphi_c, F] u_n=  \frac{c\partial_x} {\lambda-c\partial_x}(\varphi_c Fu_n)-F \frac{c\partial_x} {\lambda-c\partial_x}(\varphi_c u_n).
\end{equation} 
Thus, since $\varphi_c Fu_n\to 0$ and $\frac{c\partial_x} {\lambda-c\partial_x}\in B(L^2(\mathbb R))$ we obtain immediately that the first term of the right-hand side of (\ref{conmuta9}) goes to zero. Next, we consider
$$
w_n=\frac{c\partial_x} {\lambda-c\partial_x}(\varphi_c u_n)\equiv \mathcal P(\varphi_c u_n).
$$
We shall see that $\{w_n\}$ is bounded in $H^\alpha(\mathbb R)$. Indeed, since $\mathcal P \in B(H^\alpha(\mathbb R))\cap B(L^2(\mathbb R))$ (since $D^\alpha \mathcal P=\mathcal PD^\alpha$), we obtain 
\begin{equation}
\begin{split}
\|w_n\|_{H^\alpha} & \leqq  \|\varphi_c u_n\|_{H^\alpha}\leqq \|[D^\alpha, \varphi_c]u_n\|+ \|\varphi_c D^\alpha u_n\|+\|\varphi_c u_n\| \\
&\leqq C_0\|D^\alpha \varphi_c \|_{\infty}\|u_n\| +\|\varphi_c\|_\infty\|u_n\|_{H^\alpha}\\
&\leqq C_1\| \varphi_c\|_{H^{1+\alpha}}\| \|u_n\|_{H^\alpha}\leqq M_2,
\end{split}
\end{equation}
where we used Theorem 3.3 in \cite{Murray} and the embedding $L^\infty(\mathbb R)\hookrightarrow BMO$. Therefore, $w_n \rightharpoonup f$ in $H^\alpha(\mathbb R)$. Next, since $\mathcal P(\varphi_c u_n)\rightharpoonup 0$  in $L^2(\mathbb R)$ we obtain that $f\equiv 0$. Then, by the local compact embedding $H^\alpha(\mathbb R)\hookrightarrow L^2(\mathbb R)$ we obtain finally that $Fw_n\to 0$ as $n\to \infty$. It finishes this item and the proof of the Lemma.
\end{enumerate}
\end{proof}

\begin{obs}\label{compactp} We note that Lemma \ref{lema4} implies that the commutator operator $[\mathcal{A}^{\lambda}, F]$ is compact for $F\in C^\infty_0(\mathbb R)$.
\end{obs} 

The next lemma represents a crucial piece in the asymptotic perturbation theory.

\begin{lema}\label{lema5} Let $z\in \mathbb C$ with $Re z\leqq \frac12 c$, then there is $n>0$ such that for all $u\in C^\infty_0(|x|\geqq n)$, we have
\begin{equation} \label{Delta}
\|(\mathcal{A}^{\lambda}-z)u\|\geqq \frac14 c\|u\|,
\end{equation}
when $\lambda$ is sufficiently small.
\end{lema}

\begin{proof} Let $n>0$ such that $max_{|x|\geqq n}|\varphi_c(x)|\leqq \frac{c}{4}$. Then for $u\in C^\infty_0(|x|\geqq n)$, it follows from the proof of Lemma \ref{lema1} above  and the bounded property of the operator $\frac{c\partial_x}{\lambda+c\partial_x}$ on $L^2(\mathbb R)$, that
\begin{equation*}
\begin{array}{lll}
Re \langle (\mathcal A^\lambda-z) u, u\rangle
&\geqq (c- Re z)\|u\|^2 -Re \langle u_n\varphi_c, \frac{c\partial_x}{\lambda+c\partial_x} u_n\rangle\\
\\
& \geqq (c- Re z -\text{max}_{|x|\geqq n} |\varphi_c(x)| ) \|u\|^2\geqq \frac14 c\|u\|^2.
\end{array}
\end{equation*}
Then, since $|Re \langle (\mathcal A^\lambda-z) u, u \rangle|\leqq \|(\mathcal A^\lambda-z) u\|\|u\|$, we finish the proof.
\end{proof}

Thus, from Lemmas  \ref{lema4}-\ref{lema5} above we can apply the asymptotic perturbation theory in   \cite{VH} (see also,  Definition 19.5, Theorem 19.12 and  Lemmas 19.13-19.14 in Chapter 19 in Hislop\&Sigal \cite{hislop}) to the continuous family of closed operators $\mathcal{A}^{\lambda}$  and to get the eigenvalue perturbations of $\mathcal{A}^{0}\equiv \mathcal L_c=D^\alpha +c-\varphi_c$  to $\mathcal{A}^{\lambda}$ with a small $\lambda$ positive. More exactly, we have the following stability result for every discrete eigenvalue of  $\mathcal L_c$ with $\alpha\in (0,1)$.

\begin{teo}\label{stability} 
Each discrete eigenvalue $\gamma$ of $\mathcal L_c$ with $\gamma \leqq \frac12 c$ is stable with respect to the family $\mathcal{A}^{\lambda}$ in the sense that: there exists $\lambda_1, \delta>0$, such that for $\lambda \in (0, \lambda_1)$, we have
\begin{enumerate}
\item[(i)] $A_\delta(\gamma)=\{z: 0<|z-\gamma|< \delta\}\subset \Delta_b$, where $ \Delta_b$ is called the {\it region of boundeness} for the family $\{\mathcal{A}^{\lambda}\}$ and defined by
$$
 \Delta_b\equiv \{z: R_{\lambda}(z)= (\mathcal{A}^{\lambda}-z)^{-1}\;\;\text{exists and is uniformly bounded for}\;\; \lambda\in (0,\lambda_1)\}.
 $$

\item[(ii)] Let $\Gamma$ be a simple closed curve about $\gamma$ such that
$\Gamma\subset A_{\delta}(\gamma)\subset \rho(\mathcal{L}_c) \cap
\rho(\mathcal{A}^{\lambda})$, for all $\lambda$ small, and define
the associated Riesz projector for $\mathcal{A}^{\lambda}$
$$
P_\lambda=-\frac{1}{2\pi i} \oint_\Gamma R_{\lambda}(z) dz.
$$
Then,
\begin{equation}\label{limit5}
\displaystyle\lim_{\lambda\rightarrow0^+}||P_{\lambda}-P_{\gamma}||=0,
\end{equation}
where $P_{\gamma}$ is the Riesz projector for $\mathcal{L}_c$ and $\gamma$.
\end{enumerate}

\end{teo}

\begin{obs}\label{rem1}
It follows from Theorem  \ref{stability} that for all
$0<\lambda\ll1$, the operators $\mathcal{A}^{\lambda}$ have discrete
spectra inside the domain determined by $\Gamma$ with total
algebraic multiplicity equal to that of  $\gamma$, because  from (\ref{limit5}) we obtain that $dim(Im
P_\lambda)=dim(Im P_{\gamma})$ for $\lambda$ small. In order
to simplify the notation, we write $dim(P_\lambda)$
to refer $dim(Im P_\lambda)$.
\end{obs}

\begin{proof} As the proof follows the same lines as that of Theorem 19.12 in \cite{hislop}, by convenience of the reader we give the main points of the analysis. The prove of item (i) follows from the property of $\gamma$ to be an isolated point of the spectrum of $\mathcal L_c$ and  Lemma \ref{lema5} (see Lemma 19.14 in \cite{hislop}).

About item (ii), from  Proposition \ref{strong} and the property $\rho(\mathcal L_c)\cap \Delta_b \neq \emptyset$, we obtain from Kato \cite{kato1} the following strong resolvent convergence
$$
\lim_{\lambda \to 0^+} R_\lambda(z)u= (\mathcal L_c-z)^{-1}u, \quad \text{for all}\;\; u\in C^\infty_0(\mathbb R).
$$
Hence the Riesz projections $P_\lambda$ satisfies  
$\displaystyle\lim_{\lambda\rightarrow0^{+}}P_{\lambda}u=P_{\gamma}u$, therefore we obtain from the principle of the non-expansion of the spectrum, the inequality $\mbox{dim}(
P_{\lambda})\geq\mbox{dim}(P_{\gamma})$ (see Lemma 1.23 in  Kato
\cite{kato1}, pg. 438). Now, as $P_{\gamma}$  is self-adjoint we have $\displaystyle\lim_{\lambda\rightarrow0^{+}}P_{\lambda}^{\ast}u=P_{\gamma}u$, and so, by using Lemma 1.24 in Kato
\cite{kato1}, the two convergence of the Riesz projectors and the local compactness Lemma \ref{lema4}, we have the inequality $\mbox{dim}(P_{\lambda})\leq\mbox{dim}(P_{\gamma})$, 
for all $0<\lambda\ll1$, and so the norm convergence
of the projections in (\ref{limit5}). It finishes the Theorem.
\end{proof}

\subsubsection{The moving kernel formula}

In this subsection we study the perturbation of the eigenvalue
$\gamma=0$ associated to  $\mathcal{L}_c$  with
respect to the family $\mathcal{A}^{\lambda}$ for small $\lambda>0$.
For this purpose, we derive a moving kernel formula in the same spirit as in Lin \cite{lin} in order to
determine when the zero eigenvalue will jump to the right or to the left. We will see that it depends of the sign of the derivative of the momentum:  $\frac{d}{dc}\langle\varphi_c,\varphi_c\rangle$.

By hypotheses, we have that $\ker(\mathcal{L}_c)=\Big [\frac{d}{dx}\varphi_c\Big]$. 
 Then we obtain that  $\mbox{dim}(P_{0})=1$, for $P_0$ being the Riesz projector associated to $\gamma=0$ and $\mathcal{L}_c$.  Therefore
from Theorem \ref{stability} one has $\mbox{dim}(P_{\lambda})=1$ for
all $0<\lambda\ll1$. So, we obtain that $\mathcal{A}^{\lambda}$  has
exactly  one spectral point in the disc $B(0;\epsilon)=\{\mu\in \mathbb C:
|\mu|<\epsilon\}$,  with $\epsilon$ small, and it is non degenerate
(simple). Moreover, since the eigenvalues of $\mathcal{A}^{\lambda}$
appear in conjugates pairs, we have that there is only one real
eigenvalue $b_\lambda$ of $\mathcal{A}^{\lambda}$ inside
$B(0;\epsilon)$. We note that from the analytic property of $\mathcal{A}^{\lambda}$ with regard to $\lambda$ and from zero being a simple eigenvalue for $\mathcal L_c$, we have that the mapping $\lambda \to b_\lambda$ is analytic around zero.

The idea in the next result is to determine the sign of $b_\lambda$, for
$\lambda$ small.

\begin{lema}\label{lema6} Let $c>0$.
Assume that $\ker(\mathcal{L}_c)=\Big [\frac{d}{dx}\varphi_c \Big ]$. For
$\lambda>0$ small enough, let $b_{\lambda}\in\mathbb{R}$ be the only
eigenvalue of $\mathcal{A}^{\lambda}$ near origin. Then,
\begin{equation}\label{2eq12}
\displaystyle\lim_{\lambda\rightarrow0^+}\frac{b_{\lambda}}{\lambda^2}
= -\frac{1}{||\varphi_c'||^2}\frac{d F}{dc}
\end{equation}
with   $F(c)=\displaystyle\frac{1}{2}\langle\varphi_c,\varphi_c\rangle$. Therefore, for $\frac{dF}{dc}>0$ we obtain $b_{\lambda}<0$ and for $\frac{dF}{dc}<0$ we obtain $b_{\lambda}>0$.
\end{lema}

At this point of our theory, we will use the existence of a smooth curve of solitary wave solutions to equation $(\ref{equa})$, $c\to  \varphi_c \in H^{\alpha+1} (\mathbb R)$  (see Remark \ref{rem3} below for other version  of the limit appearing in (\ref{2eq12})). 

\begin{proof} As the proof follows the same lines as that  of  (4.7)  in Lin \cite{lin}, by convenience of the reader we give the main points of the analysis. From  Theorem $\ref{stability}$ we see
that for $\lambda>0$ small enough, there exists $u_{\lambda}\in D(\mathcal{A}^{\lambda})$,  such that 
\begin{equation}\label{norm}
\int_{\mathbb R} \varphi_c(x)|u_\lambda(x)|^2 dx=1,
\end{equation}
and 
$(\mathcal{A}^{\lambda}-b_{\lambda})u_{\lambda}=0$,
$b_{\lambda}\in\mathbb{R}$ and $\lim_{\lambda\rightarrow 0^{+}}b_{\lambda}=0$. Next, it is not difficult to see that for $Re z\leqq \frac{c}{2}$ and $u,v$ satisfying $(\mathcal{A}^{\lambda}-z)u=v$ we have
\begin{equation}\label{mainineq}
||u||_{H^{\frac{\alpha}{2}}}\leq M \int_{\mathbb R} |\varphi_c(x)||u(x)|^2 dx+ \int_{\mathbb R} |v(x)|^2 dx,
\end{equation}
with $M$ independent of $\lambda$. Thus, we obtaining immediately that $\|u_{\lambda}\|_{H^{\frac{\alpha}{2}}}\leq C$, for some constant $C>0$ which does not depend on $\lambda>0$. Therefore, $u_{\lambda} \rightharpoonup  u_0$ in $H^{\frac{\alpha}{2}}$ as  $\lambda\rightarrow0^+$. It is not difficult to see that $u_0\neq 0$ and $\int_{\mathbb R} \varphi_c(x)|u_{\lambda} -u_0|^2dx \to 0$ as $\lambda\rightarrow0^+$ (because $\varphi_c(x) \to 0$ and the renormalization condition in (\ref{norm})). Moreover, the relation $(\mathcal{A}^{\lambda}-b_{\lambda})(u_{\lambda}-u_0)=b_{\lambda}u_0 + (\mathcal L_c-\mathcal{A}^{\lambda})u_0$, implies from (\ref{norm}) and Proposition \ref{strong} that $u_{\lambda} \to u_0$ in $H^{\frac{\alpha}{2}}$. Then , $\mathcal{A}^{\lambda}u_{\lambda}\to \mathcal L_c u_0=0$ as $\lambda\rightarrow0^+$. By hypothesis we have $u_0=\theta \frac{d}{dx} \varphi_c=  \varphi'_c$ (without loss of generality we can assume $\theta =1$). So, $u_\lambda \to  \varphi'_c$ in $H^{\frac{\alpha}{2}}$. The relation,
$$
\frac{b_\lambda}{\lambda} \langle u_\lambda,  \varphi'_c \rangle=\frac{1}{c} \langle (\varphi_c -D^\alpha) u_\lambda,  \frac{c\partial_x}{\lambda-c\partial_x}\varphi_c \rangle \to -\frac{1}{c} \langle (\varphi_c -D^\alpha) \varphi'_c,   \varphi_c \rangle=  \langle  \varphi'_c, \varphi_c \rangle=0,
$$
implies that $\frac{b_\lambda}{\lambda}\to 0$ as $\lambda\to 0^+$.

Next, we calculate
$\lim_{\lambda\rightarrow0^+}\frac{b_{\lambda}}{\lambda^2}$. We write $u_\lambda=c_\lambda \varphi'_c+\lambda
v_\lambda$, with $c_\lambda=\langle u_\lambda,
\varphi'_c\rangle/\langle \varphi'_c,
\varphi'_c\rangle$. Then, $\langle v_\lambda,
\varphi'_c\rangle=0$ and $c_\lambda\to 1$ as
$\lambda\rightarrow0^+$. Following the same strategy as in \cite{lin}, we can show that $v_\lambda\to v_0$ in $H^{\frac{\alpha}{2}}$ with $v_0 \neq 0$ and satisfying $\mathcal L_c v_0= \varphi_c$ and $\langle v_0,  \varphi'_c \rangle=0$. Since  $\mathcal L_c (\frac{d}{dc} \varphi_c)=- \varphi_c$ it follows that 
$$
v_0=-\frac{d}{dc} \varphi_c+ d_0 \varphi'_c,\qquad\; d_0=\frac{\langle \frac{d}{dc} \varphi_c,  \varphi'_c \rangle}{\|\varphi'_c\|^2}.
$$
Similarly as in \cite{lin}, we can rewrite $u_\lambda=\overline{c_\lambda} \varphi'_c+\lambda
\overline{v_\lambda}$, with $\overline{c_\lambda}\to 1$, $\overline{v_\lambda}\to -\frac{d}{dc} \varphi_c$ in $H^{\frac{\alpha}{2}}$, as $\lambda\to 0^+$. Moreover, the equality
\begin{equation}\label{difi}
\frac{b_\lambda}{\lambda^2} \langle u_\lambda,  \varphi'_c \rangle=-\frac{\overline{c_\lambda}}{c}
 \langle  \frac{c\partial_x}{\lambda-c\partial_x}\varphi_c, \varphi_c \rangle-\frac{1}{c}
 \langle  \frac{c\partial_x}{\lambda-c\partial_x}(\varphi_c-D^\alpha)v_\lambda, \varphi_c \rangle
\end{equation}
implies the limit
$$
\frac{b_\lambda}{\lambda^2} \langle u_\lambda,  \varphi'_c \rangle\to \frac{1}{c}\langle \varphi_c, \varphi_c \rangle- \frac{1}{c}
 \langle (\varphi_c-D^\alpha)\frac{d}{dc} \varphi_c, \varphi_c \rangle=-  \langle \frac{d}{dc} \varphi_c, \varphi_c \rangle.
$$
Therefore,
$$
\lim_{\lambda\rightarrow0^+}\frac{b_{\lambda}}{\lambda^2}=-\frac{1}{\|\varphi'_c\|^2} \langle \frac{d}{dc} \varphi_c, \varphi_c \rangle.
$$
We finishes the Lemma.
\end{proof}

\begin{obs}\label{rem3}
In the proof of Lemma \ref{lema6} the existence of the curve of solitary waves  $c\to  \varphi_c \in H^{\alpha+1} (\mathbb R)$ was  used exactly for obtaining the  relation $\mathcal L_c (\frac{d}{dc} \varphi_c)=- \varphi_c$. Thus, it is not difficult  to see that we can change the hypothesis on the curve by the existence of $\psi\in D(\mathcal L_c)$ such that $\mathcal L_c \psi=\varphi$, with $\varphi$ being a positive solution for (\ref{equa}). Therefore, supposing that  for $
\mathcal{L_c}= D^\alpha+c -\varphi$ we have $\ker(\mathcal{L_c})=\Big [\frac{d}{dx}\varphi \Big ]$, then   relation (\ref{2eq12}) can be rewrite as
\begin{equation}\label{2eq13}
\displaystyle\lim_{\lambda\rightarrow0^+}\frac{b_{\lambda}}{\lambda^2}
= \frac{1}{||\varphi'||^2}\langle \psi, \varphi\rangle.
\end{equation}
\end{obs}

\begin{proof}$[${\bf{Theorem \ref{teo2}: Linear instability criterium  for fKdV equations}}$]$
The proof follows the same lines as in Lin \cite{lin} by using  that the mapping $\lambda\in \mathcal S \to \mathcal{A}^{\lambda}$ represents an
analytical family of operators of type-$A$, Theorem \ref{stability}, Lemma \ref{lema6}, Proposition \ref{prop1} and Lemma \ref{lema3} above. So, there exists $\lambda>0$ and $0\neq u\in H^\alpha(\mathbb R)$
such that $\mathcal{A}^{\lambda}u=0$ and therefore $e^{\lambda
t}u(x)$ is a purely growing mode solution to (\ref{lifkdv}).
\end{proof}

\begin{proof} $[${\bf{Corollary \ref{coro3}: Linear instability of ground state for fKdV equations}}$]$ From the analysis in the proof of Theorem  \ref{nonlinear}, it follows that the self-adjoint operator
$\mathcal L_c= D^\alpha +c-Q_c$ satisfies that $n(\mathcal L_c)=1$ and  $Ker(\mathcal L_c)=[\frac{d}{dx}Q_c]$. Now, the curve $c\to V(c)=Q_c\in H^\alpha (\mathbb R)$ is at least of $C^1$-class. Indeed, we know that $Q_c(x)=2c Q(c^{1/\alpha} x)$, for $Q$ being the ground state associated to (\ref{dground}), then from the relation
$$
V'(c)= 2Q(c^{1/\alpha}x) +\frac{2}{\alpha} c^{1/\alpha} xQ'(c^{1/\alpha}x)
$$
and $R=\alpha Q+xQ' \in H^{\alpha+1}(\mathbb R)$ (see proof of Theorem \ref{nonlinear}) we obtain $V'(c)\in H^{\alpha+1}(\mathbb R)\subset H^\alpha (\mathbb R)$. Thus, 
\begin{equation}
\frac{d}{dc} \langle Q_c, Q_c\rangle=4\|Q\|^2  \frac{d}{dc}c^{2-\frac{1}{\alpha}}= 4\Big(2-\frac{1}{\alpha}\Big)c^{1-\frac{1}{\alpha}} \|Q\|^2  <0
\end{equation}
where we have used that $\alpha<\frac12$ (For $\alpha=\frac12$, we obtain the equality $\frac{d}{dc} \langle Q_c, Q_c\rangle =0$ for all $c$). Hence, the condition $(ii)$ in Theorem \ref{teo2} can be applied and therefore we finish the proof.
\end{proof}

\begin{proof}$[${\bf{Theorem \ref{teo3}: Linear instability criterium  for gfKdV equations}}$]$
The proof follows the same  lines as that established for the linear instability criterium for the fKdV equation. But the strategy for showing the basic Lemma \ref{lema2} and the compactness Lemma \ref{lema4} associated to the family of linear operators $\mathcal{V}^{\lambda}:H^{\alpha}(\mathbb R)\longrightarrow L^2(\mathbb R)$, with $ \mbox{Re} \lambda
>0$, 
\begin{equation}\label{D}
\mathcal{V}^{\lambda}v\equiv cv + \frac{c\partial_x}{\lambda-c\partial_x}(\varphi_c v- \mathcal M v), 
\end{equation}
need to be changed. Indeed, with regard to Lemma \ref{lema2} for $\mathcal{A}^{\lambda}$ we change the relation (\ref{conmu2}) by
\begin{equation}\label{conmu2g}
[\mathcal V^\lambda, \zeta_d]=(1-\mathcal E^\lambda)[D^\alpha, \zeta_d] + [\mathcal E^\lambda, \zeta_d](\varphi_c-\mathcal M) + (1-\mathcal E^\lambda)[\mathcal M-D^\alpha, \zeta_d]. 
\end{equation}
The estimative required for the two first term in the right-hand side of (\ref{conmu2g}) is equal to that in (\ref{conmu3})-(\ref{conmu4}). For the third term, we use the condition in (\ref{hbeta}) for $\eta(\xi)=\beta(\xi)-|\xi |^\alpha$. Indeed, it is not  difficult  to see that the kernel
\begin{equation}\label{ker}
K_r(x,y)=(x-y)\check{\eta}(x-y)\zeta'_d(r(x-y)+y)\qquad r\in [0,1],
\end{equation}
where ``$\check{\eta}$'' represents the inverse Fourier transform of $\eta$, satisfies
\begin{equation}\label{ker2}
\int_{\mathbb R}\int_{\mathbb R}|K_r(x,y)|^2 dxdy=\int_{\mathbb R}\int_{\mathbb R}|x\check{\eta}(x)|^2|
\zeta'_d(y)|^2dxdy=\|\eta'\|^2\|\|\zeta'_d\|^2=\frac1d\|\eta'\|^2\|\|\zeta'\|^2.
\end{equation}

Therefore,
\begin{equation}\label{ker3}
\begin{array}{lll}
&\|[\mathcal M-D^\alpha, \zeta_d]u\|^2\leqq \int_{\mathbb R}|\int_0^1 \int_{\mathbb R} u(y)K_r(x,y)dy dr|^2 dx\\
\\
&\leqq \int_0^1\int_{\mathbb R}|\int_{\mathbb R} u(y)K_r(x,y)dy|^2dx dr\leqq \int_0^1 \|u\|^2 \int_{\mathbb R}\int_{\mathbb R} |K_r(x,y)|^2dydxdr\\
\\
&\leqq \frac1d\|\eta'\|^2\|\|\zeta'\|^2\|u\|^2.
\end{array}
\end{equation}
Therefore, $\|(1-\mathcal E^\lambda)[\mathcal M-D^\alpha, \zeta_d]u\|\leqq  \frac{1}{d^{1/2}}\|\eta'\|\|\|\zeta'\|\|u\|$. It finishes the estimative.

 Now, with regard to Lemma \ref{lema4} for $\mathcal{V}^{\lambda}$ we change the relation (\ref{conmuta7}) by
\begin{equation}\label{conmuta7g}
[\mathcal V^{\lambda_n}, F]=[D^{\alpha}, F] + [\mathcal W^{\lambda_n}, F] + [\mathcal M-D^{\alpha}, F].
\end{equation}
with
$$
\mathcal W^\lambda=\frac{c\partial_x}{\lambda-c\partial_x}\varphi_c-
\frac{\lambda}{\lambda-c\partial_x}\mathcal M: L^2(\mathbb R)\to L^2(\mathbb R).
$$
The estimative required for the two first term in the right-hand side of (\ref{conmuta7g}) is equal to that in the proof of Lemma \ref{lema4}. For the third term, we use the condition in (\ref{hbeta}). Indeed, since $\eta$ is bounded and continuous over $\mathbb R$ and $\eta'(x)\to 0$, as $|x|\to \infty$, (because $\eta'\in L^2(\mathbb R)$), follows from Theorem C in Cordes \cite{cords} that the commutator $[\mathcal M-D^{\alpha}, F]: L^2 (\mathbb R)\to L^2 (\mathbb R)$ is compact. It finishes the proof of the Theorem.
\end{proof}

\section{``Stability of the Blow-up" for the critical fKdV equation}

In this section we obtain information of
large-time asymptotic behaviour of solutions for the  critical fKdV equation:
\begin{equation}\label{critical}
u_t + uu_x - D^{1/2} u_x =0
\end{equation}
on $\mathbb R$.  As we saw in the last sections, the orbit
generated by ground state solutions $Q_c$ associated to equation (\ref{equa}) for $\frac12<\alpha < 2$ and $1\leqq p<p_{max}(\alpha)$ are nonlinearly stable in $H^\frac{\alpha}{2}(\mathbb R)$ by the flow of equation (\ref{0gfkdv}) for $p<2\alpha$ (see Theorem \ref{nonlinear}), moreover, the general linear instability criterium, Theorem \ref{teo2}, shows the linear instability of $Q_c$ for $\frac13<\alpha < \frac12$ and $p=1$. Now,  from the proof of these later results we can see that the  behavior of the solutions  $Q_c$ for $\alpha=\frac12$ and $p=1$ is unclear, essentially because the expression
\begin{equation}
\frac{d}{dc} \langle Q_c, Q_c\rangle=4\Big(2-\frac{1}{\alpha}\Big)c^{1-\frac{1}{\alpha}} \|Q\|^2 
\end{equation}
is zero exactly for $\alpha=\frac12$. 

Recently Saut\&Klein in \cite{KS} had provided a detailed numerical study pertaining to the dynamics of the fKdV model (\ref{0gfkdv}) with $0<\alpha<1$ and $p=1$. For the specific case of $\alpha=\frac12$ and with a initial data $u_0$ of negative energy ($E(u_0)<0$) and with a mass larger that the solitary wave mass $Q_c$ ($\|Q_c\|<\|u_0\|$) the simulations show a possible blow-up phenomenon of the associated solution (see Fig. 10 in \cite{KS}). Moreover, the peak which appears to blow-up eventually gets more and more compressed laterally, grows in hight and propagates faster with a profile of a {\it dynamically rescaled solitary wave}. Here we will show that in fact we have a kind of ``{\it stability of the blow-up}'' near to the possible unstable ground state solutions for equation (\ref{critical}).

The strategy  for showing our ``{\it stability}'' result follows that used by Angulo {\it et al.}  in \cite{ABLS} (see also Angulo \cite{angulo4}) for studying the critical case in  the model  (\ref{0gfkdv}) for $\alpha\geqq 1$, namely,  $p=2\alpha$, $p\in \mathbb N$. Thus, we consider  $\mathcal L_c$ be the linear, self-adjoint, closed,
unbounded operator defined on  $H^{1/2}(\mathbb R)$ by
\begin{equation}\label{eq4.26}
{\mathcal L_c} = D^{1/2} + c - Q_c,
\end{equation}
where $Q_c$ is the ground-state solution associated to (\ref{equa}). Therefore, from \cite{FL} we have the following properties:
\begin{enumerate}
\item[1)] ${\mathcal L_c}$ has a single negative
eigenvalue which is simple, with eigenfunction $\chi_{c}>0$, the
zero eigenvalue is simple with eigenfunction $Q_c '$, and the
remainder of the spectrum of ${\mathcal L}_c$ is positive and
bounded away from zero.
\item[2)] The curve $c\to Q_{c}$ is $C^1$ with values in
$H^{\frac32}(\mathbb R)$.
\end{enumerate}

Next, we consider the conserved energy functional $E$ in (\ref{energy}) with $\alpha=\frac12$ and $p=1$. Therefore, from (\ref{ener})-(\ref{aprior}) we have that $E(Q_c)=0$. Moreover, from (\ref{G-N}) follows the {\it a priori} estimative 
\begin{equation}
\|D^{1/4}u(t)\|^2\Big[1-\frac{\|u_0\|}{\|Q_c\|}\Big]\leqq 2E(u_0).
\end{equation}
Thus, if we consider $E(u_0)\leqq 0$ then necessarily we have the condition $\|Q_c\|\leqq \|u_0\|$

Now, we introduce the auxiliary functions
\begin{equation}\label{eq4.28}
\psi(x,t) =\mu(t)^{-\frac 12} u(\mu(t)^{-1}x, t)
\end{equation}
where
\begin{equation}\label{eq4.29}
\mu(t) = \frac{||D^{1/4}u(t)||^{4}}{||D^{1/4}
Q_c||^{4}},
\end{equation}
$\mu (0)= 1$ and $0\leqq t< t^*$ with $t^*$ the maximal time of
existence of the solution of (\ref{critical}) under consideration, if the solution is global, $t^*= +\infty$ (see Remark \ref{rnon}). Note that unless $u$
is the zero-solution, $\mu (t)\in (0,\infty)$ for $0< t< t^*$. The
normalization $\mu (0) =1$ is a temporary one made to simplify the
presentation of the argument and  it can be dispensed (see \cite{ABLS}). By
using $E$ defined in (\ref{energy}), it is easy  to check that the
function $\psi$ verifies the identities
\begin{equation}\label{propert}
\begin{array}{lll}
(i) &\|\psi(t)\|=\|u(t)\|= \|u_0\|,\\
(ii) & \langle\psi (t),D^{1/2} \psi(t)\rangle = \langle
Q_c,D^{1/2} Q_c\rangle,\\
(iii) &E(\psi(t)) =\frac 1{\mu(t)^{1/2}} E(u(t)).
\end{array}
\end{equation}

Since the stability considered here is with respect to form, i.e.,
up to translation in space, we introduce the pseudo-metric
\begin{equation}\label{eq4.33}
\begin{array}{lll} \rho_c(\psi(t),Q_c)^2 = {\inf}_{r\in \mathbb R}\,\{||D^{1/4} \psi(\cdot +r,t) -&
D^{1/4}Q_c(\cdot)||^2 +c ||\psi(\cdot+r,t)-Q_c(\cdot)||^2\}
\end{array}
\end{equation}
on $H^{1/4}(\mathbb  R)$. Define the set $\mathcal K $ to be
$$
\mathcal K=\,\{ u_0 : u_0 \in H^{s} (\mathbb R)\;\;
\text{and}\;\; E(u_0)\leqq 0\}\subset H^{\frac14} (\mathbb R),\quad s>\frac{21}{16}.
 $$
We recall that the condition $s>\frac{21}{16}$ ensures that the Cauchy problem for (\ref{critical}) is local well-posedness in $H^{s} (\mathbb R)$ (see \cite{fds2}). Of course, the problem to prove well-posedness in $H^{\frac{\alpha}{2}} (\mathbb R)$ in the general case $\alpha\in [\frac12,1)$, is still open (see Remark \ref{rnon} above).

The next theorem is a stability result which belongs to the
spatial structure of the solutions of (\ref{critical}) in the
critical case.

\begin{teo}\label{t4.9} Let $Q_{c}$ be the ground  state profile for (\ref{equa}). Then, for any $\epsilon > 0$ there is a
$\delta=\delta(\epsilon)>0$ such that if $u_0 \in \mathcal K$ with
$\rho_c (u_0,Q_c) <\delta$ and $u$ is the solution of
(\ref{critical}) corresponding to the initial value $u_0$, then $u
\in C([0, t^*); H^{1/4}(\mathbb R))$ and
\begin{equation}\label{eq4.34}
\begin{array}{lll}
{\inf}_{r\in\mathbb  R}&\Big\{c\|u(\cdot, t) -\mu(t)^{\frac12}\,
Q_c(\mu(t)(\cdot
-r))\|^2 \\
\\
\;\;\;&+\frac {1}{\mu(t)^{1/2}} \|D^{1/4} u(\cdot, t) -\mu(t)^{\frac12}
D^{1/4} Q_c(\mu(t)(\cdot -r))\|^2\Big\} <\epsilon
\end{array}
\end{equation}
for all $t \in [0,t^*)$, where $t^{*}$ is the maximal existence
time for the solution $u$  and $\mu$ is as in (\ref{eq4.29}).
\end{teo}

\begin{proof} Suppose at the outset that $\mu(0)=1$. The proof is
based on the time-dependent functional
$$
B_t[u] = \frac {1}{\mu(t)^{1/2}} E(u(t)) + \frac{c}2
\left(\frac{||u(t)||}{||Q_c||}\right)^{2k} (||u(t)||^2 -
||Q_c||^2)
$$
where $k \in \mathbb N$ will be chosen later. From the definition of
$B_t$, it is clear that if $u$ is a solution of (\ref{critical})
then $B_t[u] = B_t[u_0]$. Using (\ref{eq4.28}), (\ref{eq4.29}), and
(\ref{propert}), we may write $B_t[u]$ in terms of
$\psi$ thusly:
\begin{equation}\label{eq4.35}
\tilde B_t[\psi] = E(\psi(t)) + \frac
{c}{2}\left(\frac{||\psi(t)||}{||Q_c||}\right)^{2k}
(||\psi(t)||^2 -||Q_c||^2)
\end{equation}
where the explicit dependence on $\mu$ disappears. As it will be
argued presently, if it is established that, modulo translations,
the inequalities
\begin{equation}\label{eq4.36}
(i)\;\Delta \widetilde B_t \leqq c_0\, ||u_0-Q_c||\qquad\rm{and}
\end{equation}
\begin{equation}
(ii)\;\Delta \widetilde B_t \geqq
\,c_1\,||\psi(t)-Q_c||^2_{_\frac14} - c_2
\,||\psi(t)-Q_c||^{3}_{_{\frac14}}-\sum^{2k}_{j=1}
c_{k,j}\,||\psi(t)-Q_c||^{j+2}_{_{\frac14}},\label{eq4.37}
\end{equation}
hold for $\Delta \widetilde B_t = \widetilde B_t[\psi]-\widetilde
B_t[Q_c]$, where $c_i, c_{k,j}$ are fixed constants, then the result
in Theorem  \ref{t4.9}  follows in a well-known form. Hence, attention is turned to establishing
these bounds. The upper bound (\ref{eq4.36}) is a straightforward
consequence of $E(u_0)\leqq 0$ and $E(Q_c)=0$ (note that $\|u_0\|-\|Q_c\|\geqq 0$), where the
constant $c_0$ depends on $\|Q_c\|$ (and on an upper bound for the
choice of $\delta$). To prove (\ref{eq4.37}), consider the perturbation of
the ground state $Q_c$
\begin{equation}\label{eq4.38}
\psi(x+\gamma,t) = Q_c(x) + a(x,t),
\end{equation}
where $a$ is a real function and $\gamma=\gamma(t)$  minimizes the
functional
$$
\Pi_t(\gamma) = ||D^{1/4} \psi(\cdot + \gamma,t) - D^{1/4}
Q_c(\cdot)||^2 +c||\psi(\cdot+\gamma,t)-Q_c(\cdot)||^2.
$$
Using the representation (\ref{eq4.38}), one calculates that
\begin{equation}\label{eq4.39}
\aligned \Delta &\widetilde B_t = \widetilde B_t[Q_c+a] - \widetilde
B_t[Q_c]\\
&= E(Q_c+a) - E(Q_c) + \frac{\eta}2
\left(\frac{||Q_c+a||}{||Q_c||}\right)^{2k} (||Q_c+a||^2
-||Q_c||^2)\\
&\geqq \frac 12\langle {\mathcal L_c}a,a\rangle +
\frac{2kc}{||Q_c||^2} \langle a,Q_c\rangle^2 - c_2(c)
||a||^{3}_{_{\frac14}} - \sum^{2k}_{j=1}
c_{k,j}(c) ||a||^{j+2}_{_\frac14}.
\endaligned
\end{equation}
The inequality in (\ref{eq4.39}) is obtained using the definition
(\ref{eq4.26}) of $\mathcal L_c$, the Cauchy-Schwartz inequality and
interpolation (see (\ref{G-N})).

A suitable lower bound on the quadratic form $\mathcal L_c$ is the
next order of business. Initially, since the ground state solution $u(x,t)=Q_c(x-ct)$ is globally defined, we have from the continuous dependence theory for the model (\ref{critical}) in $H^s(\mathbb R)$, $s>\frac{21}{16}$, that for $t$ in some interval of time $[0, T]$, the inf $\Pi_t(\gamma)$ is attained in $\gamma=\gamma(t)$ for $t\in [0,T]$ (see Lemmas 6.2-6.3 in Angulo {\it et al.} \cite{ABLS}). Hence, using that $Q_c$ satisfies equation $D^{1/2} Q_c+c Q_c-\frac12 Q_c^2=0$ we obtain
$$
\frac{d}{dr} \Pi_t(r)|_{r=\gamma}=2\int_{\mathbb R} [cQ'_c(x)+D^{1/2} Q'_c(x)] a(x,t)dx=2\int_{\mathbb R} Q_c(x)Q'_c(x)a(x,t)dx,
$$
which give us the following compatibility relation on $a$, namely,
\begin{equation}\label{eq4.40}
\int_{\mathbb R} \, Q_c(x) Q_c'(x) a(x,t) \, dx =0
\end{equation}
for all $t$ in an interval of time $[0, T]$.

The issue of obtaining the lower bound (\ref{eq4.37}) for the
right-hand side of inequality (\ref{eq4.39}) is addressed in the
next few lemmas.

\begin{lema}\label{l4.10} Let ${\mathcal L_c}=D^{1/2}
+c - Q_c$. Then there exists  
$\sigma <0$ such  that if $\widetilde h =Q_c-\sigma D^{1/2} Q_c$, then
$$
{\text{min}}_{\langle f,\widetilde h\rangle =0, \|f\|=1} \langle
\mathcal L_c f,f\rangle =0.
$$
\end{lema}

{\em Proof:} For any given value $\sigma$, define the function
$f_0$ by
$$
f_0(x) = - \frac 1{c} Q_c(x) - \frac{2+2c \sigma}{c} x Q_c'(x).
$$
Then using the relation $D^{1/2}(xQ_c) =\frac12 D^{1/2} Q_c+xD^{1/2} Q_c'$, we obtain that ${\mathcal L_c} f_0 = Q_c-\sigma
D^{1/2} Q_c=\widetilde h$ and, consequently, that
$$
 \langle f_0, \mathcal L_c f_0 \rangle =\langle f_0,
\widetilde h\rangle = \left(||Q_c||^2 +\frac
1{2c}||D^{1/4} Q_c||^2\right)\sigma\\
-\frac12 ||D^{1/4}
Q_c||^2 \sigma^2.
$$
It is thus obvious that for small negative values of $\sigma$, it
is possible to have both 
$$
\langle \widetilde h, \chi_c\rangle =
\int \chi_{c}Q_c\,dx -\sigma \int \chi_c D^{1/2} Q_c \,dx \neq 0
$$
 and
\begin{equation}\label{criterio}
\langle \mathcal L_c^{-1}\widetilde h,\widetilde h\rangle=\langle
f_0,\widetilde h\rangle <0.
\end{equation}

Since $Ker(\mathcal L_c)=[Q_c']$ and $\widetilde h
\in (Ker(\mathcal L_c))^\perp$, it follows from Weinstein \cite{W1}) (see also Lemma 6.4 in \cite{angulo4}) that
\begin{equation}\label{lemwe}
\theta =\text{min}\Big\{\langle{\mathcal L_c} f,f\rangle:\,
\,||f||=1 \, \text{and} \ \langle f,\widetilde h\rangle =0
\Big\}=0.
\end{equation}
The proof of the existence of the minimum in (\ref{lemwe}) follows the same ideas as in Lemma 6.7 in Angulo {\it et al.} \cite{ABLS}.
This completes the proof of the lemma. 
\end{proof}

\begin{lema}\label{l4.11} If $\widetilde h \equiv Q_c -\sigma D^{1/2} Q_c$
with $\sigma <0$ chosen as in the last lemma, then
\begin{equation}\label{eq4.41}
\inf\Big\{\langle {\mathcal L_c}f,f\rangle :\,||f||=1, \langle f,
\widetilde h\rangle=0, f\perp Q_cQ_c'\Big\}\equiv \nu >0.
\end{equation}
\end{lema}

\begin{proof} Because of Lemma  \ref{l4.10},  it follows $\nu\geqq 0$. Suppose that $\nu=0$. Then, we can
guarantee the existence of a function $f^*$ such that $\|f^*\|=1$,
$\langle f^*,\widetilde h \rangle=0$, $ \langle f^*,Q_cQ'_c\rangle =0$
and $\langle \mathcal L_c f^*,f^*\rangle=0$. Therefore, there exists at least one
non-trivial critical point $(f^*, \tau, \theta, \nu)$ for the
Lagrange multiplier problem
\begin{equation}\label{eq4.16}
\left\{
\begin{array}{ll}
&\mathcal L_c f = \tau f +\theta \widetilde h  + \nu
Q_c Q_c',\\
&\rm{subject\; to}\\
&\|f\|=1,\,\,\langle f,Q_c Q_c',\rangle=0\;\; \rm{and}\;\; \langle f, \widetilde h\rangle=0.
\end{array}\right.
\end{equation}
Using the fact $\langle {\mathcal L}_c f^*, f^*\rangle=0$, it is
easily seen that (\ref{eq4.16}) implies $\tau =0$. Moreover,
since ${\mathcal L_c}Q_c' =0$, we have that $\langle {\mathcal L_c}
f^*, Q_c'\rangle = \langle f^*, {\mathcal L_c} Q_c' \rangle = \nu
\int (Q_c')^2 Q_cdx =0$, which implies $\nu=0$. It is
thereby concluded that
$$
{\mathcal L_c} f=\theta \widetilde h
$$
has nontrivial solutions $(f^{*},\theta)$ satisfying the 
constraints. But if $f$ is the auxiliary function arising in
the proof of Lemma  \ref{l4.10}, we have that ${\mathcal L_c}
f_0= \widetilde h$ and so ${\mathcal L_c} (f^* -\theta f_0)=0$.  Then
 $f^* -\theta f_0 \in Ker({\mathcal L_c})$. It follows from (\ref{criterio}) that
$\langle f_0, \widetilde h\rangle\neq 0$, and so $\theta =0$. Therefore,
for some non-zero $\lambda\in \mathbb R$, it is true that $f^*
=\lambda Q_c'$, which is a contradiction since such a function
cannot be orthogonal to $Q_c Q'_c$. Therefore, the minimum in
(\ref{eq4.41}) is positive and the proof of the Lemma is
completed.
\end{proof}

We note that from (\ref{eq4.41}) and from the specific form of
$\mathcal L_c$, we have that if $f\in H^{\frac14}(\mathbb R)$ satisfies
$\langle f,  \widetilde h\rangle=0$ and $\langle f,Q_cQ'_c\rangle=0$,
then
\begin{equation}\label{coerci}
\langle \mathcal L_c f,f\rangle=\int |D^{1/4}f(x)|^2 +(c-Q_c(x))|f(x)|^2dx\geqq \beta_0\|f\|^2_{\frac14},\;\;\beta_0>0.
\end{equation}

{\bf{Continuation of proof of Theorem  \ref{t4.9}}} Attention
is now turned to estimating the term $\,\frac 12\langle {\mathcal
L_c} a, a\rangle + \,\frac{2k c}{||Q_c||^2} \langle a,Q_c\rangle^2$
in (\ref{eq4.39}), where $a$ satisfies the compatibility relation
(\ref{eq4.40}). We continue to carry
 over the notation from
Lemma  \ref{l4.10}  and Lemma  \ref{l4.11}. In particular,
$\sigma$ is chosen so that the conclusions of Lemma  \ref{l4.10}
are valid. Define $a_{_{||}}$ and $a_\perp$ to be
$$
a_{_{||}} =\,\frac{\langle a,\widetilde h\rangle}{||\widetilde
h||^2}\widetilde h \;\;\;\;\rm{and} \;\;\;\;
a_\perp=a-a_{_{||}}.
$$ 
It follows from the properties of $a$ and
$\widetilde h=Q_c-\sigma D^{1/2} Q_c$ that $\langle a_\perp,
\widetilde h\rangle =0$, $\int Q_cQ_c ' a_\perp dx =0$. Without
loss of generality, take  $\langle a,\widetilde
h\rangle <0$. Thus, from Lemma  \ref{l4.11}, the Cauchy-Schwarz
inequality
 and from the properties of $a$, $a_\perp$, $a_{_{||}}$ and $\widetilde h$, it
follows that
\begin{equation}\label{eq4.42}
\left\{ \aligned \langle{\mathcal L_c}a_\perp,a_\perp\rangle &\geqq
D_1 ||a_\perp||^2,\, \, \,\langle{\mathcal
L_c}a_{_{||}},a_{_{||}}\rangle =
\,\frac{||a_{_{||}}||^2}{||\widetilde h||^2} \langle \widetilde
h,{\mathcal L_c} \widetilde h\rangle,\\
\langle{\mathcal L_c}a_{_{||}},a_\perp\rangle &= \,\frac{\langle a,
\widetilde h\rangle}{||\widetilde h||^2} \langle {\mathcal
L_c}\widetilde h, a_\perp\rangle \geqq -D_2 ||a_\perp|| \
||a_{_{||}}||
\endaligned
\right.
\end{equation}
for some positive constants $D_1$ and $D_2$. Identity $(ii)$ in 
(\ref{propert}) implies
$-2\langle a, D^\beta Q_c\rangle = ||D^{\beta/2} a||^2$.
Thus, from the Cauchy-Schwarz inequality we obtain (remember,
$\sigma$ and $\langle a,\widetilde h\rangle$ are both negative)
\begin{equation}\label{eq4.43}
\begin{array}{lll}
\frac{2kc}{||Q_c||^2} \langle a,Q_c\rangle^2
&\geqq \frac{2kc}{||Q_c||^2} \left(\langle a,\widetilde
h\rangle^2 - \sigma \langle a, \widetilde
h\rangle ||D^{1/4} a||^2\right)\\
&\geqq \frac{2kc}{||Q_c||^2} ||\widetilde h||^2 ||a_{_{||}}||^2
+2kc\sigma D_3||a||^3_{\frac14} ,
\end{array}
\end{equation}
with $D_3>0$.  We choose $\theta>0$ so that $D_1- \theta D_2
\equiv D_4
>0$. By Young's inequality, $||a_{\perp}\|||a_{_{||}}|| \leqq
\theta||a_\perp||^2 +\,\frac 1\theta||a_{_{||}}||^2$. Finally, fix
$k$ in such a way that
$$
\frac{2kc}{||Q_c||^2} ||\widetilde h||^2
+\,\frac{\langle\widetilde h, \mathcal L_c \widetilde
h\rangle}{||\widetilde h||^2} - \,\frac{D_2}{\theta} \equiv D_5
>0.
$$
With these choices, it follows from (\ref{eq4.42}) and
(\ref{eq4.43}) that
\begin{equation}\label{eq4.44}
\begin{array}{lll}
 \frac 12\langle {\mathcal L_c} a,a\rangle +
\frac{2k c}{||Q_c||^2} \langle a,Q_c\rangle^2 &\geqq D_5
||a_{_{||}}||^2 + D_4||a_\perp||^2 + 2kc \sigma
D_3||a||^3_{_{\frac14}}\\
&\geqq D' ||a||^2 - D'' ||a||^3_{_{\frac14}}
\end{array}
\end{equation}
for some positive constants $D'$ and $D''$. With (\ref{eq4.44}) in
hand, it follows easily from the specific form of the operator
$\mathcal L_c$ (see (\ref{coerci})) that
\begin{equation}\label{eq4.45}
\frac 12\langle {\mathcal L_c} a,a\rangle + \frac{2k c}{||Q_c||^2}
\langle a,Q_c\rangle^2 \geqq \widetilde D_1 ||a||^2_{_{\frac14}} -
\widetilde D_2 ||a||^3_{_{\frac14}},
\end{equation}
with $\widetilde D_1, \widetilde D_2 >0$. Finally, using
(\ref{eq4.45}) in conjunction with (\ref{eq4.39}), we obtain
$$
\aligned \Delta\widetilde B_t &\geqq \widetilde D_1
||a||^2_{\frac{1}{4}} - \widetilde D_2 ||a||^3_{_{\frac14}}-
c_2(c)||a||^{3}_{_{\frac14}}
- \sum^{2k}_{j=1} c_{k,j}(c) ||a||^{j+2}_{_{\frac14}}\\
&\geqq c_0 ||a||^2_{_{\frac14}} - c_1||a||^3_{_{\frac14}}
- \sum^{2k}_{j=1} c_{k,j} ||a||^{j+2}_{_{\frac14}}
\endaligned
$$
where $c_0,c_1,c_{k,j}$ are positive constants which depend only
on $c$.

Now we are in position to finish  Theorem  \ref{t4.9}. Suppose first
that $u_0$ lies in the set $\mathcal K$ of  ``nonpositive-energy'
initial values and suppose $||u_0 -Q_c||_{_{\frac14}} =\delta$. Then
at least for $t\in [0,T]$, it follows from (\ref{eq4.36}) and
(\ref{eq4.37}) that
\begin{equation}\label{eq4.46}
q(\rho_c(\psi(t), Q_c)) \leqq \Delta \widetilde B_t \leqq
c_0\delta
\end{equation}
where $q(x) = c_0 x^2 - c_1 x^3 -
\,\sum^{2k}_{j=1} c_{k,j} x^{j+2}$. Since
$||a(t)||^2_{_{\frac14}}$ $=\rho_c(\psi(t),Q_c)^2$ is
a continuous function of $t \in [0,t^*)$,
it follows from the inequality
\begin{equation}\label{eq4.47}
q(\rho_c(\psi(0),Q_c)) \leqq c_0\delta
\end{equation}
and (\ref{eq4.46}), that given $\epsilon >0$, then for all $t \in
[0, T]$,
\begin{equation}\label{eq4.48}
\rho_c(\psi(t), Q_c) \leqq \epsilon,
\end{equation}
provided that  $\delta$ is chosen small enough at the outset. To
finish the proof, we need to show that inequality (\ref{eq4.48})
is still true for $t \in [0,t^*)$. This part is shown using a method similar
to that of the proof of Theorem  6.1 in Angulo {\it et al.} \cite{ABLS}. Therefore, the
stability in Theorem  \ref{t4.9} is established if $\mu(0)=1$. The
general case, wherein the initial data is not necessarily such that
$\mu(0)=1$  requires a little more of work, and therefore we refer the reader to see the
reference \cite{ABLS}. This completes the proof of Theorem
\ref{t4.9}.

\subsection{Behaviour of the stability parameters for the critical-fKdV equation} 

In the proof of Theorem \ref{t4.9}, we use that there is a specific choice of the  translation parameter $\gamma=\gamma(t)$ such that 
\begin{equation}\label{beha}
||D^{1/4} \psi(\cdot +\gamma,t) -
D^{1/4}Q_c(\cdot)||^2 +c ||\psi(\cdot+\gamma,t)-Q_c(\cdot)||^2 =\rho_c(\psi(t),Q_c)^2\leqq \epsilon
\end{equation}
for all $t<t^*$, where $\psi$ is the rescaled version of the solution $u$ of (\ref{critical}) defined in (\ref{eq4.28}). Moreover, a choice of $\gamma$ for which (\ref{beha}) holds may be determined via the orthogonality condition in (\ref{eq4.40}). By an application of the implicit-function theorem as in  Lemma 4.2  in \cite{ABLS}, it is obtained that as long as $\psi$ satisfies (\ref{beha}), {\it there is a unique, continuously differentiable choice of the value $\gamma(t)$ that achieves (\ref{eq4.40})} provided that the initial data $u_0\in H^s(\mathbb R)$ for $s$ sufficiently large and the profile $Q_c\in H^n(\mathbb R)$ for $n$ large (at least for $n\geqq 3$). Moreover, with the hypothesis of sufficiently regularity for the initial data $u_0$ we can see that $\mu$ defined in (\ref{eq4.29}) belongs to the class $C^1([0, t^*): \mathbb R)$.

Thus, by following the line of argumentation in Lemma 4.3 in \cite{ABLS}, we obtain the relation between the translation and dilation parameters involved in our stability result in Theorem \ref{t4.9}.

\begin{teo}\label{beha2} Let $Q_{c}$ be the ground  state profile for (\ref{critical}) such that $Q_c\in H^n(\mathbb R)$, $n\geqq 3$. Then, for any $\epsilon > 0$ there is a
$\delta=\delta(\epsilon)>0$ such that if $u_0 \in H^s(\mathbb R)\cap \mathcal K$, with $s$ sufficiently large 
and $\|u_0-Q_c\|_{\frac14}<\delta$, then there exists a $C^1$-mapping $\gamma: [0, t^*)\to \mathbb R$ such that
\begin{enumerate}
\item[i)] $\|\psi(\cdot+\gamma(t), t)-Q_c\|_{\frac14}\leqq \epsilon$\quad for $t\in [0,t^*)$,
\item[ii)] for all $t\in [0,t^*)$,
$$
\Big | \gamma(t)-c \mu(t)\int_0^t \sqrt{\mu(s)}ds\Big|\leqq C \epsilon \mu(t)\Big(\int_0^t \sqrt{\mu(s)}ds +\int_0^t \frac{|\mu'(s)|}{\mu^2(s)}ds\Big)
$$
where $C$ depends only on $Q_c$.
\end{enumerate}
\end{teo}

\begin{obs}\label{beha3} The statement in Theorem \ref{beha2} deserves to be clarified at  least in some points and its relation with the nonlinear stability result established in Theorem \ref{nonlinear} above for the case $\alpha\in (\frac12, 2)$, $p< p_{max}(\alpha)$ and $p<2\alpha$.
\begin{enumerate}
\item [(1)] The regularity required on the initial data $u_0$ is given to ensure that the associated solution $u$ satisfies in a classical sense the equation (\ref{critical}).

\item [(2)] An similar analysis may be made for obtaining the behavior of the parameter of translation involved  in the nonlinear stability result in Theorem \ref{nonlinear}. In this case, $\mu(t)\equiv 1$ for all $t$ and so for $\gamma=\gamma(t)$ such that 
\begin{equation}\label{beha4}
\|u(\cdot+\gamma, t)-Q_c\|_{\frac{\alpha}{2}}\leqq \epsilon,
\end{equation}
satisfies for all $t\in [0,t^*)$,
$$
\Big | \gamma(t)-ct\Big|\leqq C \epsilon t
$$
where $C$ depends only on $Q_c$.
\end{enumerate}
\end{obs}

\section{Nonlinear stability and linear instability for the fBBM equation} 

This section is devoted to the fractional BBM equation
\begin{equation}\label{fBBM}
u_t+u_x + \partial_x(u^2)+D^\alpha u_t=0,
\end{equation}
for $\alpha\in (\frac13, 1)$. As the structure of the analysis is similar to that used for the fKdV, we will only indicate the new basic differences.

Consider a solitary wave solution $u(x,t)= \psi_c(x-ct)$, $c>1$, of the fBBM equation(\ref{fBBM}). Then the profile $\psi_c$ satisfies the equation
\begin{equation}\label{sbbm}
D^\alpha \psi_c+ \Big(1-\frac1c\Big) \psi_c -\frac1c \psi^2_c=0.
\end{equation}
Therefore, we obtain the following Pohozaev identity
\begin{equation}\label{pohobbm}
(3\alpha-1)\int_{\mathbb R} |D^{\alpha/2} \psi_c|^2 dx=\Big(1-\frac1c\Big)\int_{\mathbb R} |\psi_c|^2 dx,
\end{equation}
proving that no finite energy solitary waves exist when $c>1$ and $\alpha\leqq \frac13$ hold. 

By considering the new variable
$v(x,t)=u(x+ct,t)-\psi_c(x)$, it follows from
(\ref{fBBM}) that
\begin{equation}\label{fBBM2}
(\partial_t-c\partial_x)(v+D^\alpha v)+\partial_x(v+2\psi_c v+O(\|v\|^2))=0.
\end{equation}
The equation
\begin{equation}\label{linearfBBM1}
(\partial_t-c\partial_x)(v+D^\alpha v)+\partial_x(2\psi_c v+v)=0,
\end{equation}
represents the linearized equation for (\ref{fBBM}) around of
$\psi_c$. So, we will give sufficient conditions
for obtaining that the solution $v\equiv 0$ is unstable by the
linear flow of (\ref{linearfBBM1}). More exactly, we are interested
to find a \textit{growing mode solution} of (\ref{linearfBBM1}) with
the form $v(x,t)=e^{\lambda t}u(x)$ and $\mbox{Re}\lambda>0.$
Thus, we obtain that $u$ satisfies the following non-local differential equation,
\begin{equation}\label{modeeqgfBM1}
D^\alpha u+u+\displaystyle\frac{\partial_x}{\lambda-c\partial_x}(u+2\psi_cu)=0.
\end{equation}
This motivates us to define the following family of closed linear operators
$\mathcal{B}^{\lambda}:H^\alpha(\mathbb R)\longrightarrow L^2(\mathbb R)$, $ \mbox{Re} \lambda>0$, given by
\begin{equation}\label{modeoperfBBM}
\mathcal{B}^{\lambda}v\equiv (D^\alpha+1)v+\displaystyle
\frac{\partial_x}{\lambda-c\partial_x}(v+2\psi_cv).
\end{equation}
Next, we consider the unbounded self-adjoint operator $\mathcal L_0: H^\alpha(\mathbb R)\longrightarrow L^2(\mathbb R)$ associated to (\ref{sbbm})
\begin{equation}\label{selfbbm}
\mathcal L_0=D^\alpha + \Big(1-\frac1c\Big)  -\frac2c \psi_c,
\end{equation}
and so $\psi'_c\in Ker(\mathcal L_0)$.

Our first result is about the behavior of $\mathcal B^\lambda$ by depending of $\lambda$.

\begin{prop}\label{strong2} For $\lambda>0$, the operator $\mathcal B^\lambda$ converges to $\mathcal L_0$ strongly in $L^2(\mathbb R)$ when $\lambda \to 0^+$, and  converges to $D^\alpha + 1$ strongly in $L^2(\mathbb R)$ when $\lambda \to \infty$.
\end{prop}

\begin{proof} Similar to that of Proposition \ref{strong}.
\end{proof} 

Next, we   localized  the essential spectrum of $\mathcal B^\lambda$, $\sigma_{ess}(\mathcal B^\lambda)$. 

\begin{prop}\label{prop2fbbm}
For any $\lambda>0$, we have
\begin{equation}
\sigma_{ess}(\mathcal B^\lambda)\subset \big \{ z: Re z\geqq \frac12 \Big(1-\frac1c\Big )\Big\}.
\end{equation}
\end{prop}

The idea of the proof of Proposition \ref{prop1} is the same of Proposition 1 in Lin \cite{lin}. The next lemma is similar to Lemma \ref{lema2} above.

\begin{lema}\label{lema2fbbm}
Given $\lambda>0$. Let $\zeta\in C_0^\infty(\mathbb R)$ be a cut-off function such that $\zeta|_{\{|x|\leqq R_0\}} =1$, for some $R_0$. Define $\zeta_d(x)=\zeta(x/d)$, $d>0$. Then, for each $d$, the operator $\zeta_d(\mathcal B^\lambda-z)^{-1}$ is compact for some $z\in \rho(\mathcal B^\lambda)$, and  there exists $C(d)\to 0$ as $d\to \infty$ such that for any $u\in C_0^\infty(\mathbb R)$,
\begin{equation}\label{conmuta}
\|[\mathcal B^\lambda, \zeta_d]u\|\leqq C(d) (\|\mathcal B^\lambda u\|+\|u\|).
\end{equation}
\end{lema}

Next, we study the behavior of $\mathcal B^\lambda$ near infinity. The next result shows the non-existence of growing modes at the left-hand  side of the complex-plane for large $\lambda$ (see Lin \cite{lin}), so, since the eigenvalues of $\mathcal B^\lambda$ appear in conjugate pairs, there are not  growing modes in all for large $\lambda$.

\begin{lema}\label{lema3fbbm} There exists $\Lambda>0$, such that for $\lambda>\Lambda$, $\mathcal B^\lambda$ has no eigenvalues in $\{z: Re z\leqq 0\}$.
\end{lema}

 Next, we study the behavior of $\mathcal B^\lambda$ for small positive $\lambda$.  It result extends those of Lin in \cite{lin} for the case of the fBBM equation (\ref{fBBM}).

\begin{lema}\label{lema4fbbm} Given $F\in C_0^\infty (\mathbb R)$. Consider any sequence $\lambda_n\to 0^+$ and $\{u_n\} \subset  H^\alpha(\mathbb R)$ satisfying
\begin{equation}\label{conmuta4}
\|\mathcal B^{\lambda_n} u_n\|+\|u_n\|\leqq M_1<\infty
\end{equation}
for some constant $M_1$. Then if $w-\lim_{n\to \infty} u_n=0$, we have
\begin{equation}\label{conmuta5}
\lim_{n\to \infty}\|Fu_n\|=0
\end{equation}
and
\begin{equation}\label{conmuta6}
\lim_{n\to \infty} \|[\mathcal B^{\lambda_n}, F] u_n\|=0.
\end{equation}
\end{lema}

\begin{proof} The proof in this case is immediate. First, $[\mathcal B^{\lambda}, F]= [D^\alpha, F] +[\mathcal G^\lambda, F]$ with
$$
\mathcal G^\lambda= \frac{\partial_x} {\lambda-c\partial_x}(1+2\psi_c).
$$
Now, from the proof of Lemma \ref{lema4} above we have that $[D^\alpha, F]$ is a compact operator. Next, the convergence
$$
[\mathcal G^{\lambda_n}, F]u_n\to 0,\qquad \text{as}\;\; n\to \infty,
$$
is proved by following the same ideas in proof of Lemma 2.5 in \cite{lin}. It finishes the Lemma.
\end{proof}

The next lemma is basic in the  stability theory (see also  Chapter 19 in Hislop\&Sigal \cite{hislop}), and its proof follows from the estimative
\begin{equation*}
Re \langle (\mathcal B^\lambda-z) u, u\rangle\geqq \frac14 \Big(1-\frac1c\Big)\|u\|^2.
\end{equation*}

\begin{lema}\label{lema5fbbm} Let $z\in \mathbb C$ with $Re z\leqq \frac12 (1-\frac1c)$, then there is $n>0$ such that for all $u\in C^\infty_0(|x|\geqq n)$, we have
\begin{equation} \label{Deltafbbm}
\|(\mathcal{B}^{\lambda}-z)u\|\geqq \frac14 c\|u\|,
\end{equation}
when $\lambda$ is sufficiently small.
\end{lema}

Thus, from Lemmas  \ref{lema4fbbm}-\ref{lema5fbbm} above, we  have the following stability result for every discrete eigenvalue of  $\mathcal L_0$ with $\alpha\in (0,1)$.

\begin{teo}\label{stabilityfbbm} 
Each discrete eigenvalue $\kappa_0$ of $\mathcal L_0$ with $\kappa_0 \leqq \frac12(1-\frac1c)$ is stable with respect to the family $\mathcal{B}^{\lambda}$. 
\end{teo}

Next, we establish the moving kernel formula for the fBBM equations.

\begin{lema}\label{lema6fbbm} Let $c>1$.
Assume that $\ker(\mathcal{L}_0)=\Big [\frac{d}{dx}\psi_c \Big ]$. For
$\lambda>0$ small enough, let $\kappa_{\lambda}\in\mathbb{R}$ be the only
eigenvalue of $\mathcal{B}^{\lambda}$ near origin. Then,
\begin{equation}\label{2eq12fbbm}
\displaystyle\lim_{\lambda\rightarrow0^+}\frac{\kappa_{\lambda}}{\lambda^2}
= -\frac{1}{c}\frac{1}{||\psi_c'||^2}\frac{d M}{dc}
\end{equation}
with   $M(c)=\displaystyle\frac{1}{2}\langle (D^\alpha + 1)\psi_c,\psi_c\rangle$. Therefore, for $\frac{dM}{dc}>0$ we obtain $\kappa_{\lambda}<0$ and for $\frac{dM}{dc}<0$ we obtain $\kappa_{\lambda}>0$.
\end{lema}

\begin{proof} The proof is similar to that of Lemma 2.7  in Lin \cite{lin}. 
\end{proof}

The linearized instability result for the  fBBM equation (\ref{fBBM})
is the following: 

\begin{teo}\label{lfbbm}$[${\bf{Linear instability criterium for fBBM equations}}$]$
Let $c\to  \psi_c \in H^{\alpha+1} (\mathbb R)$ be a smooth curve of positive solitary wave solution to
equation $(\ref{sbbm})$ with $\alpha \in (\frac13, \frac12)$, $p=1$. The wave-speed $c$ can be considered over some nonempty interval $I$, $I\subset (1, +\infty)$. We assume that the self-adjoint operator $\mathcal L_0= D^\alpha + \Big(1-\frac1c \Big)-\frac2c\psi_c$ with domain $D(\mathcal L_0)= H^ \alpha(\mathbb R)$ satisfies
\begin{equation}\label{eq1400}
Ker(\mathcal{L}_0)=[\frac{d}{dx}\psi_c].
\end{equation}
Denote by $n(\mathcal{L}_0)$ the number (counting multiplicity)
of negative eigenvalues of the operator $\mathcal{L}_0$.  Then there is a purely growing mode
$e^{\lambda t}u(x)$ with $\lambda>0$, $u\in H^s(\mathbb R) -\{0\}$, $s\geqq 0$, to the linearized equation $(\ref{linearfBBM1})$ if one
of the following two conditions is true:
\begin{enumerate}
\item[(i)] $n(\mathcal{L}_0)$ is even and $\frac{d}{dc}M(c)>0$.
\item[(ii)] $n(\mathcal{L}_0)$ is odd and $\frac{d}{dc}M(c)<0$,
\end{enumerate}
where $M(c)=\displaystyle\frac{1}{2}\langle (D^\alpha + 1)\psi_c,\psi_c\rangle$. 
\end{teo}

Next, we show the  nonlinear stability and linear instability theorems for the fBBM established in the introduction (Theorems \ref{nonlinearfbbm1}-\ref{linsfbbm1}) for the ground state solutions associated to the  equation (\ref{sbbm}). Before, we study the sign of the quantity $\frac{d}{dc}M(c)$.

\begin{lema}\label{derfbbm} Let $c>1$ and $\alpha\in (\frac13, 1)$.
 Then for any solution $\psi_c$ of (\ref{sbbm}) we have 
\begin{equation}\label{M} 
2M(c)=\Big[ \frac{1}{3\alpha-1}c^{\frac1\alpha-1}(c-1)^{3-\frac1\alpha} +c^{\frac1\alpha}(c-1)^{2-\frac1\alpha}\Big] \|\Psi\|,
\end{equation}
where $\Psi$ satisfies $D^\alpha \Psi +\Psi-\Psi^2=0$. Thus, we obtain that
\begin{equation}\label{dM}
\frac{d}{dc}M(c)= \left\{
\begin{array}{lll}
>0\qquad \text{for}\;\; \alpha\in [\frac12, 1)\;\;\text{and}\;\; c>1\\
\\
>0\qquad \text{for}\;\; \alpha\in (\frac13, \frac12)\;\;\text{and}\;\; c>c_0\\
\\
<0\qquad \text{for}\;\; \alpha\in (\frac13, \frac12)\;\;\text{and}\;\; 1<c<c_0.
\end{array}\right.
\end{equation}
where $c_0>1$ is the bigger positive root of the polynomial $q(c)=6\alpha^2 c^2-4 c\alpha +1-\alpha$, and it is given in (\ref{c_0}) below.

\end{lema}

\begin{proof} Initially, suppose that $\psi_c$ is a solution of (\ref{sbbm}) with $c>1$ then for the scaling $\Psi(x)= a\psi_c(bx)$ with
\begin{equation}\label{ab}
a=\frac{1}{c-1},\qquad \text{and}\qquad b=\Big (\frac{c}{c-1}\Big)^{1/\alpha}
\end{equation}
we obtain that the profile $\Psi$ satisfies $D^\alpha \Psi+\Psi-\Psi^2=0$. Thus, we obtain  that
$$
(D^\alpha +1) \psi_c(x)= \frac{1}{ab^\alpha} D^\alpha \Psi(x/b)+\frac{1}{a}\Psi(x/b).
$$
Therefore, from the relation $(3\alpha-1)\int |D^{\alpha/2} \Psi(x)|^2dx=\int |\Psi(x)|^2dx$ we obtain
\begin{equation}\label{eM}
2M(c)= \Big[\frac{1}{3\alpha-1} c^{\frac{1}{\alpha}-1}(c -1)^{3-\frac{1}{\alpha}} +c^{\frac{1}{\alpha}}(c -1)^{2-\frac{1}{\alpha}}\Big ] \|\Psi\|^2\equiv p(c) \|\Psi\|^2.
\end{equation}
Next, we determine the sign of the derivative of $p(c)$ defined in (\ref{eM}). An simple calculation shows that
\begin{equation}\label{eM2}
p'(c)=\frac{c^{\frac{1}{\alpha}}(c-1)}{(c-1)^{\frac{1}{\alpha}}}\Big[ \frac{1-\alpha}{\alpha}\frac{1}{3\alpha-1}\frac{(c-1)^2}{c^2} +\frac{2}{\alpha}\frac{c-1}{c} +\frac{2\alpha-1}{\alpha}\Big ].
\end{equation}
Thus, it follows immediate from (\ref{eM2}) that for $c>1$ and $\alpha\in [\frac12, 1)$ we have $p'(c)>0$. Next, it is no difficult to see that $ p''(c)>0$ for every $c>1$. Moreover, since $\alpha<\frac12$ and 
$$
\lim_{c\to 1^{+}} \frac{c^{\frac{1}{\alpha}}(c-1)}{(c-1)^{\frac{1}{\alpha}}}=+\infty
$$
we have $\lim_{c\to 1^{+}} p'(c)=-\infty$. Now, for $c\to +\infty$
$$
p'(c)\approx c\Big [\frac{1-\alpha}{\alpha}\frac{1}{3\alpha-1} + \frac{2\alpha+1}{\alpha}\Big ]
$$
and for $1>\alpha>\frac13$, we get that $\lim_{c\to +\infty} p'(c)=+\infty$. Therefore, there is an unique point $c_0>1$  such that $p'(c_0)=0$. Thus, we obtain for $\alpha\in (\frac13, \frac12)$ that $p'(c)<0$ for $c\in (1,c_0)$ and $p'(c)>0$ for $c\in (c_0, +\infty)$. 

Now, for determining $c_0$, we have that $p'(c)=0$ if and only if $q(c)=6\alpha^2 c^2-4 c\alpha +1-\alpha=0$. Since, $q(1)=(3\alpha-1)(2\alpha-1)<0$ we have that the real zeros  of $q$, $r_0$ and $c_0$, satisfy $r_0<1<c_0$. The  exact value of $c_0$ is given by
\begin{equation}\label{c_0}
c_0=\frac{2+\sqrt{2(3\alpha-1)}}{6\alpha}.
\end{equation}
It finishes the proof.
 \end{proof}

\begin{proof} $[${\bf{Proof of Theorem \ref{nonlinearfbbm1}}}$]$ For $\alpha \in (\frac13, 1)$, the scaling $Q(x)=a\Phi_c(bx)$, with $a$ and $b$ defined in (\ref{ab}), implies that the ground state $Q$ satisfies $D^\alpha Q+Q-Q^2=0$. Therefore, from \cite{FL} follows that the self-adjoint operator $\mathcal L= D^\alpha +1-2Q$
satisfies $Ker(\mathcal L)=[\frac{d}{dx}Q]$ and $n(\mathcal L)=1$. Thus, a similar analysis as that in the proof of Theorem \ref{nonlinear} above we conclude that for $\mathcal J_c= D^\alpha +\Big (1-\frac1c \Big )-\frac2c \Phi_c$ satisfies $Ker(\mathcal J_c)=[\frac{d}{dx}\Phi_c]$ and $n(\mathcal J_c)=1$.

Now, for  $M(c)=\displaystyle\frac{1}{2}\langle (D^\alpha + 1)\Psi_c,\Psi_c\rangle$ we have from Lemma \ref{derfbbm} that $M'(c)>0$ exactly for $\alpha\in [\frac12, 1)$ and $c>1$, and for $\alpha\in (\frac13, \frac12)$ and $c>c_0$. Hence, from Grillakis {\it et.al} theory \cite{grillakis1} we finish the proof.
\end{proof}

As consequence of Theorem \ref{lfbbm} we obtain the prove of Theorem \ref{linsfbbm1}. 

\begin{proof}$[${\bf{Proof of Theorem \ref{linsfbbm1}}}$]$  From the proof of Theorem \ref{nonlinearfbbm1} we have that the self-adjoint operator $\mathcal J_c= D^\alpha +\Big (1-\frac1c \Big )-\frac2c \Phi_c$ satisfies $Ker(\mathcal J_c)=[\frac{d}{dx}\Phi_c]$ and $n(\mathcal J_c)=1$. Moreover, Lemma \ref{derfbbm} establishes that $\frac{d}{dc} M(c)<0$ for $\alpha\in (\frac13, \frac12)$ and $c\in (1, c_0)$, where $c_0$ is giving in (\ref{c_0}). It finishes the Theorem.
\end{proof}

\begin{obs}\label{refbbm} Next we have the following observations about the stability Theorems \ref{nonlinearfbbm1}-\ref{linsfbbm1}.
\begin{enumerate}
\item [(1)]  Similarly as in the case of the fKdV, the statement of stability in Theorem \ref{nonlinearfbbm1} is a {\bf{conditional}} one, in  the sense that for all $\epsilon>0$ there is a $\delta>0$ such that if $u_0\in  H^s(\mathbb R) \cap U_\delta $, for $s>\frac32 -\alpha$, then $u(t)\in U_\epsilon$, for all $t\in (-T_s, T_s)$, where $T_s$ is the maximal time of existence of $u$ satisfying $u(0)=u_0$. We recall that the best known result of local well-posedness for the fBBM model (\ref{fBBM}) is given in \cite{fds2} for initial data in $H^s(\mathbb R)$, $s>\frac32 -\alpha$ and $\alpha\in (0,1)$. It which does not allow to globalize the solution using conservation laws.

\item [(2)]  We recall that in \cite{fds} was showed the existence and stability of solitary waves solutions for the fBBM  by considering the minimization problem
$$
I_q=inf \Big\{\int_{\mathbb R} u^2 +|D^{\alpha/2} u|^2 dx: u\in H^\frac{\alpha}{2}(\mathbb R)\;\;\text{and}\;\; \int_{\mathbb R} \frac{u^2}{2} + \frac{u^3}{3}dx=q\Big \}.
$$
For $\alpha \in (\frac13, \frac12)$, a critical value constrain  $q_0=q_0(\alpha)$ was established in such way that for $q >q_0$, the set of ground state solutions associated to the variational problem above will be stable in $H^\frac{\alpha}{2}(\mathbb R)$. From our analysis in Lemma \ref{derfbbm}  and Theorem \ref{nonlinearfbbm1}, we note that that critical value constrain  $q_0$ can be determined explicitly in terms of the threshold value $c_0$ in (\ref{c_0}).

\item[(3)] Our orbital stability  and linear instability results  in Theorem \ref{nonlinearfbbm1} and Theorem \ref{linsfbbm1} for the fBBM equation show a scenario similar  to that known for the generalized BBM equation (GBBM)
$$
u_t+u_x + u^pu_x-u_{xxt}=0.
$$
Indeed, the critical exponent for the stability of solitary waves solutions for the GBMM is $p=4$, though the explanation for instability when $p\geqq 4$ is different. In fact, from Souganidis\&Strauss \cite{SS} solitary waves of the GBBM of arbitrary positive velocity are stable when $p<4$ but when $p\geqq 4$ there exists $c_\ast=c_\ast(p)$ such that the solitary waves of velocity $c<c_\ast$ are unstable (nonlinearly) while those of velocity $c>c_\ast$ are nonlinear stable.

\end{enumerate}
\end{obs}

\vskip0.2in

\indent\textbf{Acknowledgements:}  This work was done while the author was visiting the Department of Mathematics of Paris-Sud University as a visiting professor and him was supported by FAPESP (S\~ao Paulo Research Fundation/Brazil) under the process 2016/07311-0. The author would like express their  thanks to the Professors Felipe Linares and  Jean-Claude Saut for useful discussions and comments  in the development of the work.

\end{document}